\documentclass[a4paper, 12pt]{article}
\usepackage{enumerate, theorem}
\usepackage{amsmath, amsfonts, amssymb}
\usepackage[height=22.5cm, width=15cm]{geometry}
\usepackage[all]{xy}
\usepackage{mathrsfs, yfonts}
\usepackage[utf8]{inputenc}

\newcommand{\Abc}{\tilde{\mathcal{A}}_{conv.}}

\newcommand{\C}{\mathbb{C}}

\newcommand{\E}{\mathcal{E}}
\newcommand{\F}{\mathcal{F}}

\newcommand{\h}{\mathcal{H}}

\newcommand{\parag}[1]{\paragraph{\sc{#1.}}}

\newtheorem{thm}{Theorem}[subsection]
\newtheorem{defn}[thm]{Definition}
\newtheorem{cor}[thm]{Corollary}
\newtheorem{prop}[thm]{Proposition}
\newtheorem{lemma}[thm]{Lemma}

\begin{document}
\title{Generalized Brieskorn Modules II:\\
Higher Bernstein Polynomials and Multiple Poles of $\frac{1}{\Gamma(\lambda)} \vert f\vert^{2\lambda}$}

\author{Daniel Barlet\footnote{Barlet Daniel, Institut Elie Cartan UMR 7502  \newline
Universit\'e de Lorraine, CNRS, INRIA  et  Institut Universitaire de France, \newline
BP 239 - F - 54506 Vandoeuvre-l\`es-Nancy Cedex.France. \newline
e-mail : daniel.barlet@univ-lorraine.fr}.}

\date{29 F\'evrier 2025}

\maketitle

\hfill {\it   Si vous ne r\'eussisez pas  \`{a} tout int\'egrer} \\
$\qquad \qquad.\qquad \qquad  \qquad\qquad \qquad \qquad\qquad \qquad  \quad$ {\it  essayer donc d'int\'egrer donc par parties}.

\parag{Abstract} 

Our main result is to show that the existence of  a root in $-\alpha - \mathbb{N}$ for the $p$-th Bernstein polynomial of the (a,b)-module generated by a holomorphic form $\omega \in \Omega^{n+1}_0$ in the (convergent)  Brieskorn (a,b)-module $H^{n+1}_0 $ associated to $f$, under the hypothesis that $f$ has an isolated singularity at the origin  relative to the eigenvalue $\exp(2i\pi\alpha)$ of the monodromy, produces  poles of order at least $p$ for the meromorphic extension of the (conjugate) analytic functional, for some $h \in \mathbb{Z}$:
$$ \omega' \in \Omega^{n+1}_0 \mapsto \frac{1}{\Gamma(\lambda)}\int_{\mathbb{C}^{n+1}} \vert f\vert^{2\lambda} \bar f^{-h} \rho\omega\wedge \bar \omega' $$
at points $-\alpha - N$ for $N$ and $h$ well chosen integers. This result is new, even for $p = 1$.
As a corollary, this implies that, in the case of an isolated singularity for $f$,  the existence of a root in $-\alpha -\mathbb{N}$ for the $p$-th Bernstein polynomial of the (a,b)-module generated by a holomorphic form $\omega \in \Omega^{n+1}_0$ implies the existence of at least $p$ roots (counting multiplicities) for the usual reduced Bernstein polynomial of the germ $(f, 0)$.\\
In  the case of an isolated singularity for $f$, we obtain that for each $\alpha \in ]0, 1] \cap \mathbb{Q}$ the biggest root $-\alpha -m$ of the reduced Bernstein polynomial of $f$ in $-\alpha - \mathbb{N}$ produces a pole  at $-\alpha -m$ for some $h \in \mathbb{Z}$ for the meromorphic extension of the  distribution 
$$\square \longrightarrow  \frac{1}{\Gamma(\lambda)}\vert f\vert^{2\lambda} \bar f^{-h}\square.$$

\parag{AMS classification}  32 S 25; 32 S 40 ; 34 E 05

\parag{Key words} Generalized Brieskorn Modules;  Asymptotic Expansions ; Period-integral;  Convergent (a,b)-Module; Hermitian Period;  Geometric (convergent) (a,b)-Modules; (convergent) Fresco; Bernstein Polynomial ; Higher  Bernstein Polynomials. \\

\tableofcontents

\newpage 

\section{Introduction}

\subsection{The aim of this article}

The roots of the reduced  Bernstein polynomial $b_{f,0}$ of the germ of holomorphic function at the origin in $\mathbb{C}^{n+1}$ control the poles of the meromorphic extension of the distribution 
$$ \square \longrightarrow \frac{1}{\Gamma(\lambda)}\int_{\mathbb{C}^{n+1}} \vert f\vert^{2\lambda} \bar f^{-h} \square $$
which is defined in a neighborhood of $0 \in \mathbb{C}^{n+1}$ (see for instance \cite{[B.81]} or \cite{[Bj.93]}).\\
The first goal of this article is to show that, assuming that $0$ is an isolated singularity for the eigenvalue $\exp(2i\pi\alpha)$ of the monodromy (this corresponds to our hypothesis $H(\alpha,1)$), the roots of the Bernstein polynomial of the (a,b)-module generated by the germ $\omega$ of holomorphic $(n+1)$-form at the origin in the Brieskorn (a,b)-module $H^{n+1}_0$  of $f$ at $0$, control the poles of the (conjugate) analytic functional defined  on $\Omega^{n+1}_0 $ by polar parts of  poles in $-\alpha - \mathbb{N}$ of the meromorphic functions
$$ \omega' \mapsto \frac{1}{\Gamma(\lambda)}\int_X \vert f\vert^{2\lambda}\bar f^{-h}\rho\omega\wedge \bar \omega' $$
where $\rho \in \mathscr{C}^\infty_c(\mathbb{C}^{n+1})$ is identically $1$ near $0$ and with a sufficiently small support (note that the polar parts  of these meromorphic extensions at points in $-\alpha - \mathbb{N}$ are independent of the choices of $\rho$ thanks to our hypothesis $H(\alpha,1)$).\\

Our goal is to give a sufficient condition, still on the (a,b)-module generated by the germ $\omega$, to obtain higher order poles for such integrals.\\
The difficulty commes now from the fact that it is not clear when, for instance,  two roots, $-\alpha-m$ and $-\alpha - m'$ with $m,m' \in \mathbb{N}$, of the Bernstein polynomial give a simple pole or a double pole for such a meromorphic extension at points $-\alpha - N$, for some choice of $\omega'$ and for some  integers $N$ and $h$  well chosen.\\
So we try  to understand when such a pair of roots are ``linked'', so produces a double pole for some choice of $\omega', h$ and $N$, or are  ``independent'', so produces at most a simple pole for any choices of $\omega', h$ and $N$. \\
As it is known that the nilpotent part of the monodromy is related to this phenomenon (see \cite{[B.84-a]} and \cite{[B.84-b]}) we have defined in our previous paper 
\cite{[part I]}, the action of the monodromy on  a simple pole geometric\footnote{The (a,b)-modules deduced from the Gauss-Manin connection which appear here are always geometric. See below} (a,b)-modules and we have shown in {\it loc. cit.}  that the natural semi-simple filtration of a geometric (a,b)-module $\mathcal{E}$ is  related to the filtration induced by the nilpotent part of the action of the monodromy on its saturation $\mathcal{E}^\sharp$ by $b^{-1}a$.\\
The section 2 of this article is devoted to some reminder of part I (see \cite{[part I]}) and to the definition and the study of the {\bf higher  Bernstein polynomials} of a geometric (a,b)-module. The case of frescos which is used in our main result is examined in details in section 3.

\subsection{The main results}

We shall describe now our main results, which are obtained  in section 4, using the tools  introduced in our previous paper \cite{[part I]} and the first sections.\\
The proofs will be given in section 4 and the reader will find the notations and definitions used in thie following presentation  explained further in this article.\\

We consider a germ $f$ of holomorphic function at the origin of $\mathbb{C}^{n+1}$ with an isolated singularity at $0$  for the eigenvalue $\exp(2i\pi \alpha)$ of its monodromy. That is to say with the hypothesis $H(\alpha, 1)$ where $\alpha$ is in $]0, 1] \cap \mathbb{Q}$.\\
We denote $f : X \to D$ a Milnor's representative of the germ $f$ near $0$.\\
In this situation we consider the {\bf generalized Brieskorn module} $H^{n+1}_0$ which is the quotient of the $(n+1)$-cohomology group of the complex $(Ker^\bullet df_0, d)$ by its  $b$-torsion.
This complex is the sub-complex of the germ at $0$ of the holomorphic de Rham complex at $0$ where $Ker^p df_0$ is the kernel of $\wedge df : \Omega^p_0 \to \Omega^{p+1}_0$ for $p \geq 2$ and $Ker^1df $ is the quotient of the kernel of the map  $\wedge df : \Omega^1_0 \to \Omega^2_0$ by $\mathbb{C}df$. Its cohomology groups after quotient by their torsion are geometric $B[a]$-modules  (see \cite{[B.I]} or  \cite{[B.II]}).\\ 
For $\omega \in \Omega^{n+1}_0$ we denote by $\F_\omega$ the fresco generated by $[\omega]$ in $H^{n+1}_0$; it is defined as 
 $\F_\omega = B[a]\omega \subset H^{n+1}_0$ where $B := \mathbb{C}\{\{b\}\}$ and $ab - ba = b^2$. \\
We refer the reader to section 2 for the definitions of the Bernstein polynomial $B_\E$ and of the definition of  {\bf higher Bernstein polynomials} $B^p_\E$ of a geometric (a,b)-module $\E$. Remind that $B_\E$ always divides the product of the $B^j_\E$. Remind also that any root of a $B^j_\E$ is a root of $B_\E$ and that  $B_\E$ is the product of the $B^j_\E$ for all $j$ when $\E$ is a {\bf fresco}. The nilpotent order $d(\E)$  of a geometric (a,b)-module $\E$ is the smallest integer $d\geq 0$ such that $B^{d+1}_\E \equiv 1$.\\

Remind that for a geometric (a,b)-module $\E$ there exists a maximal quotient $\E^{[\alpha]}$ of $E$ which is $[\alpha]$-primitive. This means that $B_{\E^{[\alpha]}}$ is the quotient of $B_\E$ by all its roots which are not in $-\alpha -\mathbb{N}$. Then  $B^p_{\E^{[\alpha]}} $ is also  the quotient of $B^p_\E$ by all its roots which are not in $-\alpha -\mathbb{N}$ for each $p$.\\

\begin{thm}\label{1}
Under the hypothesis $H(\alpha, 1)$ for $f$, let $\omega \in \Omega^{n+1}_0$ such that  the nilpotent order of the fresco $\F^{[\alpha]}_\omega$  is equal to $p \in \mathbb{N}^*$, where  $\F_\omega := B[a]\omega \subset H^{n+1}_0$. Then there exists $\omega' \in \Omega^{n+1}_0$ and  $h \in \mathbb{N}$ be such that the meromorphic extension of the function 
\begin{equation}
F^{\omega, \omega'}_h(\lambda) := \frac{1}{\Gamma(\lambda)}\int_X \vert f \vert^{2\lambda} \bar f^{-h} \rho \omega\wedge \bar \omega' 
\end{equation}
holomorphic for $\Re(\lambda) \gg 1$, has a pole of order at least $p$ at a point in $-\alpha - \mathbb{N}$, where $\rho \in \mathscr{C}^\infty_c(X)$ is identically $1$ near $0$.\\
Conversely, if the the meromorphic extension of $F^{\omega, \omega'}_h(\lambda) $, for some choice of $\omega'$ and $h$, has a pole of order $p$ at a point in $-\alpha - \mathbb{N}$ then the nilpotent order of the fresco  $\F^{[\alpha]}_\omega$ (and then also the nilpotent order of  $\h^{[\alpha]}$) is at least equal to $p$.
\end{thm}

The following Lemma (see \cite{[part I]}) allows to  enlighten how the previous result give a relation between the nilpotent order of $\h^{[\alpha]}$ and the poles  in $-\alpha - \mathbb{N}$ of the meromorphic extension of the distributions $\frac{1}{\Gamma(\lambda}\vert f\vert^{2\lambda} \bar f^{-h}$ under our hypothesis $H(\alpha, 1)$.

\begin{lemma}\label{1 bis}
For any geometric (a,b)-module $\E$ such that $-\beta$ is a root of $B^p_\E$ there exists $x \in \E$ and an integer  $q \geq p$ such that $B^q_{\F_x}$ has a root  
$-\beta + m$ with $m \in \mathbb{N}$, where $\F_x := B[a]x \subset E$ is the fresco generated by $x$ in $\E$.\\
Moreover, if $p = d(\E^{[\alpha]})$ where $[\alpha]$ is the class of $\beta$ in $\mathbb{Q}\big/\mathbb{Z}$ we may choose $x$ to have $m = 0$.
\end{lemma}

\begin{cor}\label{1 ter}
In the situation of the theorem above, is $p$ is the nilpotent order of $\h^{[\alpha]}$ there exists $\omega \in \Omega^{n+1}_0$ which satisfies the hypothesis of the theorem above.
 Moreover, if $-\beta$ is a root of $B^p_{\h^{[\alpha]}}$ we may fing $\omega \in \Omega^{n+1}_0$ such that $-\beta$ is a root of $B^p_{\F_\omega}$.
\end{cor}

The next result precises that the biggest root in $-\alpha - \mathbb{N}$ of $B^p_{\F_\omega}$ is a pole for the function $F^{\omega, \omega'}_h(\lambda)$ for a convenient choice of $\omega'$ and $h$.

\begin{thm}\label{2}
Under the hypothesis $H(\alpha, 1)$ for $f$,  let $\omega \in \Omega^{n+1}_0$ such that  the nilpotent order of $\F^{[\alpha]}_\omega$, where  $\F_\omega := B[a]\omega \subset H^{n+1}_0$,  is equal to $p \in \mathbb{N}^*$. Let $-\alpha -m$ be the biggest root of $B^p_{\F_\omega}$ inside $-\alpha - \mathbb{N}$. Then there exists $\omega'$ and $h$ such that $F^{\omega, \omega'}_h(\lambda)$ has a pole of order $p$ at $\lambda = - \alpha -m$.
\end{thm}

The last result (see Corollary \ref{improve 2} below)  strengthens strongly  the main theorem  in \cite{[B.22]}.

\begin{thm}\label{3}
Under the hypothesis $H(\alpha, 1)$ for $f$, let $\omega \in \Omega^{n+1}_0$ be such that the nilpotent order of $\F^{[\alpha]}_\omega$ is equal to $p$.  For each $s \in [1, p]$  let $\xi_s$  be  the biggest number
 in $-\alpha - \mathbb{N}$ for which there exists $\omega'$ and $h$ such that $F^{\omega, \omega'}_h(\lambda)$ has a pole of order $\geq s$ at $\xi_s$. Then there exists $j \in \mathbb{N}$ such that $\xi_s$ is a root of $B^{s+j}_{\F_\omega}$.\\
 Moreover, if $\xi_s = \xi_{s+1} = \dots = \xi_{s+h}$, then there exist at least $h$  distinct integers $q_1, \dots, q_h $  such that $\xi_s$ is a root of $B^{s+q_j}_{\F_\omega}$ for $j \in [1, h]$.\\
 \end{thm}
 
 The proofs  and some more precise results are given in Sections 4.4 and 4.5.

\newpage

 \section{Asymptotics and geometric (a,b)-modules}

Let $B := \C\{\{b\}\}$ be the algebra of constant coefficients microdifferential operators of degree $\leq 0$. So $b := \partial_s^{-1}$ and a series
$$ S(b) := \sum_{p=0}^\infty  c_pb^p \quad {\rm is \ in} \  B \ {\rm when} \quad \exists R > 1  \  \exists C_R  \quad {\rm such \ that} \ \forall p \geq 0\quad  \vert c_p\vert \leq C_R R^p p!  .$$
Then $B$ acts on $\C\{s\}$, the algebra of germs of holomorphic functions at the origin in $\C$, by the rule
$$ S(b)[f] = \sum_{n=0}^\infty \frac{1}{n!}\Big(\sum_{p = 0}^n  (n-p)! c_p\gamma_{n-p} \Big)s^n \quad {\rm where} \ f(s) := \sum_{q=0}^\infty \gamma_q s^q .$$
Note that the derivation in the variable $b$ acts continuously on $B$\footnote{For its natural dual Fréchet topology associated to the  "pseudo-norms" :\\
$\vert\vert \quad \vert\vert_R := Sup  \{ \vert c_p\vert R^{-p}/p! , \ p \geq 0\}.$} and that the multiplication by $s$ on $\C\{s\}$, denoted by $a$, satisfies the commutation relation
$$ aS(b) = S(b)a + b^2S'(b) \qquad \forall S \in B$$
in the algebra of continuous endomorphisms of $\C\{s\}$ denoted by $End_c(\C\{s\})$.
We shall denote by $B[a]$ the sub-algebra of  $End_c(\C\{s\})$ which is generated by $B$ and $a := \times s$.\\
We shall denote by $A := \C\{a\}$ the sub-algebra of $End_c(\C\{s\})$ of elements given by multiplication by some $f \in \C\{s\}$.

Define 
$$ e_{\alpha, j} := s^{\alpha - 1}\frac{(Log\, s)^j}{j!} $$

 for $\alpha \in ]0, 1]  \cap \mathbb{Q}$ and for $j \in \mathbb{N}$ and let $\Xi_\alpha^{(N)}$ be the free $\C\{s\}$-module with basis $e_{\alpha, j}$ for $j \in [0, N]$. Then define the action of $B$ on $\Xi_\alpha^{(N)}$ inductively  by the formulas
 $$  \alpha be_{\alpha, j} =  s e_{\alpha, j} -  be_{\alpha, j-1} \qquad {\rm for} \quad  j \geq 1 \qquad {\rm and} \quad  \alpha be_{\alpha, 0} = s e_{\alpha, 0} $$
 with the commutation relations
 $$ S(b)a = aS(b)  - b^2S'(b)  \qquad \forall S \in B$$
 where $S'$ is the derivative (in the variable $b$) in the algebra $B$.\\
 
 \begin{lemma}\label{26/1}
 For each $\alpha  \in ]0, 1]  \cap \mathbb{Q}$ and each $N \in \mathbb{N}$ the action of $\C[b]$ on  $\Xi_\alpha^{(N)}$ extends to an action of $B$ and makes  it isomorphic to  the free $B$-module with basis $e_{\alpha, j}$ for $j \in [0, N]$, on which the action of $a$ is continuous and satisfies $ab - ba = b^2$.\\
 So it is a left $B[a]$-module for which the action of $\C[a]$ extends continuously  to an action of $\C\{a\}$.
 \end{lemma}
 
 For a proof, see section 2.3 in \cite{[part I]}.$\hfill \blacksquare$\\
 

 \begin{defn}\label{27/1}
 For $\alpha \in ]0, 1[\cap \mathbb{Q}$ we define
 $$ S\Xi_\alpha^{(N)} := \Xi_\alpha^{(N)} $$
 and for $\alpha = 1$
 $$ S\Xi_1^{(N)} := \Xi_1^{(N+1)}\big/\Xi_1^{(0)} .$$
 Then for $\mathscr{A}$ a finite subset in $]0, 1]\cap \mathbb{Q}$, $N \in \mathbb{N}$ and $V$ a finite dimensional complex vector space, define $S\Xi_\mathscr{A}^{(N)}\otimes V$ as before, replacing $\Xi_\alpha^{(N)}$ by $S\Xi_\alpha^{(N)}$ for each $\alpha \in \mathscr{A}$.
 \end{defn}
 Note that $\Xi_1^{(0)} = \C\{s\}$;  so it  is the non singular part (i.e. the uni-valued  holomorphic part) of the expansion at the origin.\\

From now on a {\bf  geometric (a,b)-module}, we also use the terminology {\bf Generalized Brieskorn Modules}, will be, by definition, a sub $B[a]$-module of some $S\Xi^{(N)}_R \otimes V$ where $\mathscr{A}$ is a finite subset in $]0, 1] \cap \mathbb{Q}$, $N$ a non negative integer and $V$ a finite dimensional vector space.\\

 Of course, thanks to the results in \cite{[part I]} section 5, this is equivalent to the "standard definition" given in {\it loc. cit.} of a
convergent geometric  (a,b)-module.

\parag{Remarks}
\begin{enumerate}
\item  Since $B$ is local noetherian and $S\Xi^{(N)}_\mathscr{A} \otimes V$ is free and finite type over $B$, any geometric (a,b)-module is a free finite type module over $B$ with a continuous $\mathbb{C}$-linear action of  $a$.
\item Also any $B[a]$-sub-module of a geometric (a,b)-module is again a geometric (a,b)-module.
\item Fix $\alpha \in ]0, 1]$.  For each integer $n$ we have in $\Xi_\alpha^{(0)}$:
$$b^ns^{\alpha-1} = \frac{s^{\alpha + n - 1}}{\prod_{p=1}^n (\alpha + p - 1)}$$
and so
$$ \Big(\sum_{n=0}^\infty c_n b^n\Big)s^{\alpha-1} = \Big( \sum_{n=0}^\infty \frac{c_n}{\prod_{p=1}^n (\alpha + n -1)} a^n\Big) s^{\alpha-1} $$
and the estimates
$$ \frac{\vert c_n\vert}{\prod_{p=2}^n (\alpha+p-1)} R^n (n-1)!  \leq \vert c_n\vert R^n \leq   \frac{\vert c_n\vert}{\prod_{p=1}^n (\alpha+p-1)} R^n n!$$
shows that the topology of $\Xi_\alpha^{(0)}$ as a free rank $1$ $A$-module or as free rank $1$ $B$-module are the same.\\
This extends easily to any $S\Xi_\mathscr{A}^{(N)} \otimes V$ and shows that in  any sub-$B[a]$-module of a $S\Xi_\mathscr{A}^{(N)} \otimes V$ the action of $a$ extends continuously to
an action of $A = \mathbb{C}\{a\}$. So any geometric  convergent (a,b)-module is a free  finite type $A$-module.

 \end{enumerate}
 
 \subsection{Some facts from \cite{[part I]}}

We shall use the following definitions  which are  equivalent to the definitions given and used in {\it loc. cit.} where $E$ is  a geometric (a,b)-module. The reader may find proofs there.
 \begin{enumerate}
\item A sub-module $F \subset E$ is {\bf normal} when $F \cap bE = bF$. In this case the quotient $E/F$ is again without $B$-torsion so free and finite type as a $B$-module and it is again a geometric (a,b)-module\footnote{This point is not obvious with the definition adopted  here but is rather easy with the standard definition.}.
\item For an arbitrary sub-module $F$ in $E$ we denote by 
$$N_E(F) = \{x \in E\ \exists n \in \mathbb{N}  \ such \ that \  b^nx \in F \}.$$
 We call this sub-module the {\bf normalisation} of $F$ in $E$. It is the smallest normal sub-module in $E$ which contains $F$.

\item  We denote by  $E^\sharp$  the {\bf saturation} of $E$ in $b^{-1}a$. \\
Note that any $S\Xi^{(N)}_R \otimes V$ is stable by $b^{-1}a$. So $E^\sharp$ is again a finite type free $B$-module and there exists an integer $m \geq 0$ such that 
$b^mE^\sharp \subset E$. So $E^\sharp\big/E$ is a finite dimensional $\C$-vector space.\\
We shall say that $E$ has a {\bf simple pole} when $E = E^\sharp$. This is equivalent to the fact that $aE \subset bE$ (since $a$ is injective).\\
  \item The rank $1$ geometric (a,b)-modules are classified by positive rational numbers. To $\alpha \in \mathbb{Q}^{*+}$ corresponds 
  $E_\alpha := B[a]\big/B[a](a-\alpha b)$. Then $E_\alpha$ is the $B$-module $Be_\alpha$ where the action of $a$ is defined by $aSe_\alpha = \alpha bSe_\alpha + b^2S'e_\alpha$ for $S \in B$, where $S'$ is the derivative (in $b$) of $S$.
 \item  For a given $\alpha \in ]0, 1] \cap \mathbb{Q}$ we denote by $E^{[\alpha]}$ the quotient of $E$ by elements in $E$ that present no non zero term like 
 $s^{\alpha+m-1}(Log\, s)^j\otimes v$ for any integers $j$ and $m$ and any $v \in V$. Since we make the quotient  of $E$  by $E \cap (S\Xi^{(N)}_{\mathscr{A}\setminus \alpha} \otimes V)$ and since  we have 
 $$ S\Xi^{(N)}_\mathscr{A} \otimes V \big/ S\Xi^{(N)}_{\mathscr{A}\setminus \alpha} \otimes V \simeq S\Xi^{(N)}_{\alpha} \otimes V$$
 the quotient $E^{[\alpha]}$ is again a geometric (a,b)-module.  A geometric  (a,b)-module $E$ such that $E = E^{[\alpha]} $  will be called {\bf $[\alpha]$-primitive} and, for a  general $E$,  $E^{[\alpha]}$ will be called {\bf the $[\alpha]$-primitive quotient} of $E$ (see \cite{[part I]} Proposition 3.3.3).\\
 When $E$ has a simple pole, $E$ is the direct sum of the $E^{[\alpha]}$ when $\alpha$ describe the image in $ \mathbb{Q}\big/\mathbb{Z} \simeq  ]0, 1] \cap \mathbb{Q}$ of the roots of $B_E$ (see below).
 \item The (usual)  {\bf Bernstein polynomial}  $B_E$  of the geometric (a,b)-module $E$ is, by definition the minimal polynomial of the endomorphism $-b^{-1}a$ acting on the finite dimensional  complex vector space
 $E^\sharp\big/bE^\sharp$ whose dimension is equal to the $B$-rank of $E$.
  \parag{Remark} If we have an exact sequence of {\bf simple poles} (a,b)-modules
 $$0 \to E_1 \to E \to E_2 \to 0 $$
 then $B_{E_1}$ and $B_{E_2}$ divide $B_E$ since the maps in  the exact sequence of finite dimensional vector spaces 
 $$ 0 \to E_1/bE_1 \to E/bE \to E_2/bE_2 \to 0$$
 commute with the respective  actions  $-b^{-1}a$ on these spaces.
 
\parag{Warning} Without the assumption that $E$ has a simple pole (which implies that $E_1$ and $E_2$ have also  simple poles) it is still true that $B_{E_2}$ divides $B_E$ since the induced map $E^\sharp \to E_2^\sharp$ is still surjective. \\
 But the kernel of this induced map is not equal to $E_1^\sharp$ in general (it is equal to the normalisation of $E_1^\sharp$ in $E^\sharp$),  so
 it is not true in general that $ B_{E_1}$ divides $B_E$. See, for instance, the case where $E$ is a fresco described below.

 \item A  {\bf fresco}\footnote{See \cite{[B.09]} for a detailed study of frescos.} is a geometric (a,b)-module $F$  such that the (finite dimensional)  complex vector space $F\big/ (aF + bF)$ has dimension $1$ (note that this vector space  is a quotient of $F/bF$). This is equivalent to ask that $F$ is a geometric (a,b)-module which  is isomorphic to the quotient of $B[a]$ by a left ideal. \\
  For a fresco the Bernstein polynomial is equal to  the {\bf characteristic polynomial} of $-b^{-1}a$ acting on $F^\sharp\big/bF^\sharp$. So its degree is the $B$-rank of $F$.
\item  A normal sub-module $G$ of a fresco $F$ is  a fresco and the quotient $F/G$ is again a fresco. 
 \item For each element $x$  in a geometric (a,b)-module $ E \subset S\Xi_\mathscr{A}^{(N)} \otimes V$ we define $d(x)$ as the maximal integer such that a non zero  term like  $s^{\alpha+m-1}(Log\, s)^{j-1}\otimes v$ appears in $x$ for some $\alpha \not= 1$ or like
  $s^{m}(Log\, s)^{j} \otimes v$  for $\alpha = 1$ for any $m \in \mathbb{N}$ and any $v \in V$.\\
   The integer $d(x)$ is called the {\bf nilpotent order} of $x$ and does not depend of the realization of $E$ as a sub-module of some $S\Xi_R^{(N)} \otimes V$.\\
   Then we define $d(E) := \sup_{x \in E} \{d(x)\}$ and call this integer {\bf the nilpotent order} of $E$.
   \item The following properties for a geometric (a,b)-module $E$ are equivalent:
   \begin{enumerate}[(i)]
   \item $d(E) \leq 1 $.
   \item $E$ is a sub-module of a finite direct sum $\oplus_{j=1}^p E_{\alpha_j}$.
   \item $E^\sharp$ is isomorphic to a finite direct sum $\oplus_{j=1}^p E_{\alpha_j}$.
   \end{enumerate}
   A geometric (a,b)-module is {\bf semi-simple} when it satisfies one of the properties above.\\
  Note that the Bernstein polynomial of a semi-simple geometric (a,b)-module has only simple roots (obvious from property $(iii)$ above).
  \item We define the {\bf semi-simple filtration}\footnote{See \cite{[part I]} Section 4  for more details.} $\big(S_j(E), j \in \mathbb{N}\big)$ by
  $$ S_j(E) := \{ x \in E \ / \ d(x) \leq j\} .$$
  Then $S_0(E) = \{0\}, S_{d(E)}(E) = E$ and for $j \in [0, d(E)-1]$ we have a strict inclusion  $S_j(E) \subsetneq S_{j+1}(E)$. Moreover, for each integer $j$, $S_j(E)$ is a normal sub-module of $E$.\\
  A geometric (a,b)-module is semi-simple when $d(E) \leq 1$. For instance, each quotient $S_j(E)\big/S_{j-1}(E)$ is semi-simple because  another way to define the semi-simple filtration is  to see that $S_j(E)\big/S_{j-1}(E)$ is the semi-simple part of $E\big/S_{j-1}(E)$ that is to say the subset of $x$ in this quotient with $d(x) \leq 1$.
 \item For any sub-module $G$ of a geometric (a,b)-module $E$ and any integer $j$ we have $S_j(G) = S_j(E) \cap G$.\\
 When $G$ is normal, the induced map $S_j(E) \to S_j(E/G)$ is not surjective in general when $j < d(E)$. For instance, it will be proved below that for $G = S_1(E)$ the image by the map induced by the  quotient map  $S_j(E) \to S_j(E/G)$ is equal to $S_{j-1}(E/G)$ (see Lemma \ref{21/1} below.)
 
 \item   It is proved in \cite{[part I]} Proposition 4.2.8   that the rank (as $B$-module) of $S_{j+1}(E)\big/S_j(E)$ is non increasing in $j$.\\
  \end{enumerate}
   
   The following easy lemma will be used later on.
   
   \begin{lemma}\label{21/1}
   Let $E$ be a geometric (a,b)-module and put $G := E/S_1(E)$. Then for each positive integer there exists a natural isomorphism of $B[a]$-modules 
   $$S_{h+1}(E)/S_1(E) \to S_h(G) $$ for each $h \geq 0$.  This implies the isomorphisms of $B[a]$-modules
     $$ S_{h+1}(E)\big/S_h(E) \to S_h(G)\big/S_{h-1}(G) .$$
 for each $h \geq 1$.$\hfill \blacksquare$\\
   \end{lemma}
   
   \parag{Proof} It is clear  with the definition of a geometric (a,b)-module adopted here that if $x$ belongs to $S_{h+1}(E)$ then its image in $G$ is in $S_h(G)$. This gives an injective morphism
   $$ S_{h+1}(E)\big/S_1(E) \to S_h(G) $$
   which is clearly surjective. $\hfill \blacksquare$\\
 
 \subsection{Higher Bernstein Polynomials}
 
 \begin{defn}\label{0}
 Let $E$ be a geometric (a,b)-module. For any $j \in [1, d(E)]$ we denote $B^j_E$ the Bernstein polynomial of the semi-simple (a,b)-module $S_j(E^\sharp)\big/S_{j-1}(E^\sharp)$.
 We call it {\bf the $j$-th Bernstein polynomial of $E$}.
 \end{defn}
 
 \parag{Remarks}\begin{enumerate}[(i)]
 \item Like $B_E$ the (usual) Bernstein polynomial of $E$, the polynomials $B^j_E$ only depend on $E^\sharp$ the saturation of $E$ by $b^{-1}a$.
 \item Since for each $j \in [1, d(E)]$ the geometric (a,b)-module $S_j(E^\sharp)\big/S_{j-1}(E^\sharp)$ is semi-simple and has a simple pole, it is a direct sum of rank $1$ geometric (a,b)-modules. So each $B^j_E$ has only simple roots.
 \end{enumerate}

 The first interesting point of this definition is given by the following proposition..
 
\begin{prop}\label{2}  Let E be a  geometric  (a,b)-module. Then for each integer $ j $ in $[1,d(E)]$ the polynomial $B^j_E$ 
 divides $B_E$. Moreover each root of $B_E$ is a root of $B^j_E$ for at least one $j \in [1,d(E)]$.
\end{prop}

\parag{Proof} We have an exact sequence of simple poles (a,b)-modules
 $$ 0 \to S_j(E^\sharp)\big/S_{j-1}(E^\sharp) \to E^\sharp\big/S_{j-1}(E^\sharp) \to E^\sharp\big/S_j(E^\sharp)\ \to 0 $$
 which gives the fact that a root of $B^j_E$ is a root of the quotient $ E^\sharp\big/S_{j-1}(E^\sharp)$ and then a root of $B_E =  B_{E^\sharp}$.\\

We have also  an exact sequence of simple pole (a,b)-modules
$$ 0 \to S_{j-1}(E^\sharp)\to S_j(E^\sharp) \to S_j(E^\sharp)/S_{j-1}(E^\sharp) \to 0$$
so it is clear that $B^j_E$  divides $B_{S_j (E^\sharp)}$.  
This already proves that for $j= d(E) = d(E^\sharp)$ the polynomial $B^d_E$  divides $B_E$. Applying this to $ F := S_j(E^\sharp)$  gives, since $d(F) = j$ and
$$S_{j-1}(F) = S_{j-1}(E^\sharp) \cap F = S_{j-1}(E^\sharp) $$
 that $B^j_E$ divides $B_F$. But $F$ is a normal sub-module of $E^\sharp$ so, thanks to point 6 above ($F$ has a simple pole) 
  $B_F$ divides $B_E$.$\hfill \blacksquare$\\

   \begin{prop}\label{22/12}
  Let $E$ be a geometric (a,b)-module and put $d := d(E)\geq 2$. If $-\beta$ is a root of $B^d_E$ there exists at least  a  root  of $B^{d-1}_E$ in $-\beta + \mathbb{N}$.
  \end{prop}
  
  \parag{Proof} It is of course enough to consider the case where $E$ has a simple pole by definition of the standard and higher Bernstein polynomials.\\
  We are only concerned by the $[\beta]$-primitive part of $E$ so we may assume that $E$ is $[\beta]$-primitive, and we may assume that $-\beta$ is the biggest root of $B^d_E$ in $-\beta + \mathbb{Z}$.\\
  
  Since $E/S_{d-1}(E)$ is a direct sum of $E_{\alpha_j}$,  we may always solve  in $E/S_{d-1}(E)$ the equation $ax - (\beta-p)bx = by$ for any $p \in \mathbb{N}^*$ and any 
  $y \in E/S_{d-1}(E)$ because $b^{-1}a - (\beta - p)$ is invertible in $E/S_{d-1}(E)$ for each integer 
   $p \geq 1$\footnote{To show that in $E_\alpha := Be_\alpha$ we can always solve the equation $(a - (\alpha-p)b)x = by$ for each given $y \in E_\alpha$ reduces to show that for each $S \in B$ we may find $T \in B$ such that $bT'(b) + pT(b) = S(b)$ for each $p \in \mathbb{N}^*$. This is an easy exercise.}.\\
  But our hypothesis implies the existence of $x_0 \in E \setminus (bE + S_{d-1}(E))$ such that 
  $$ax_0 - \beta bx_0 = b^2y_0 + z_0 $$
  for some $y_0 \in E$ and $z_0 \in S_{d-1}(E)$,  which is equivalent  to 
  $$(b^{-1}a - \beta)(x_0) \in bE + S_{d-1}(E) $$  
   since the simple pole of $E$ implies that we may write $z_0 = bz_1$ with $z_1 \in S_{d-1}(E)$, since  $z_0$ is in  $bE \cap S_{d-1}(E) = bS_{d-1}(E)$, as $S_{d-1}(E)$ is normal in $E$. \\
   Now there exists $y_1 \in E$ such that $ay_1 - (\beta-1)by_1 = by_0 + \xi_0$ with $\xi \in S_{d-1}(E)$ (inversiblity of $b^{-1}a - (\beta - 1)$ in $E\big/S_{d-1}(E)$). And again we may write $\xi_0 = b\xi_1$ with $\xi_1 \in S_{d-1}(E)$. Applying $b$ to this equality gives
   $$ aby_1 - \beta b^2y_1 = b^2y_0 + b^2\xi_1 $$
   and then
   $$  a(x_0 - by_1) - \beta b(x_0 - by_1) = bz_1 - b^2\xi_1 \in bS_{d-1}(E).$$

   \smallskip
   
   From now on, assume that $d = 2$.\\
   If we assume that there is no root of the Bernstein polynomial of $S_{d-1}(E) = S_1(E)$ in $-\beta + \mathbb{N}$ we may solve in $S_1(E)$, which is semi-simple and has a simple pole (so it is  a direct sum of $E_{\alpha_j}$),  the equation
   $$ az_2 - \beta bz_2 = b(z_1 - b\xi_1)$$
   with $z_2 \in S_1(E)$. Then we find that $x = x_0 - by_1 - z_2$ satisfies $ax - \beta bx = 0$. Since $x$ is not in $bE$, $Bx$ is a normal rank $1$ sub-module of $E$ isomorphic to $E_\beta$, so it is contained in $S_1(E) $. Contradiction !\\
    So there exists a root in $-\beta + \mathbb{N}$ for $S_{d-1}(E)$. Since we assume  $d = 2$ the proof is complete for this case. \\
    
    Now for $d(E) \geq 3$ we have the equalities, where we define  $F := E/S_{d-2}(E)$:
    $$ S_1(F) = S_{d-1}(E)/S_{d-2}(E) , \quad d(F) = 2 \quad {\rm and} \quad B^d_E =  B^2_F , \quad B^{d-1}_E = B^1_F$$
    and applying to $F$ the previous result gives that $B^{d-1}_E$ has a root in $-\beta + \mathbb{N}$.$\hfill \blacksquare$\\
    
    \begin{cor}\label{23/12}
    Let $E$ be a geometric (a,b)-module and let $j$ be an integer in $[1, d(E)]$. Then if $-\beta$ is a root of $B^j_E$,  for each $h \in [1, j]$ there exists a root of 
    $B^h_E$ in $-\beta + \mathbb{N}$.
    \end{cor}
    
    \parag{Proof} It is enough to apply the previous proposition successively  to each simple pole   (a,b)-module $S_j(E^\sharp)$ since the equality $S_{j-h}(S_j(E^\sharp)) = S_{j-h}(E^\sharp)$ for each $h \in [0,j-1]$  shows the equality  $B^{j-h}_{S_j(E^\sharp)} = B^{j-h}_E$. $\hfill \blacksquare$\\
    
    Note that if the $j$ roots find in $-\beta + \mathbb{N}$ are two by two distinct, we have found $j$ roots of $B_E$, since each $B^h_E$ divides $B_E$.
    
    
     \begin{cor}\label{29/1-b}
 Let $E$ be a geometric (a,b)-module and assume that $-\beta$ is the greatest root of $B_E$ in $-\beta + \mathbb{Z}$. Then $-\beta$ is a root of $B^1_E$.\\
 Moreover, if $-\beta$ is a root of $B^j_E$ then $-\beta$ is a root of $B^h_E$ for each $h \in [1, j]$.
 \end{cor}
 
 \parag{Proof} This is an obvious consequence  of Corollary \ref{23/12}.$\hfill \blacksquare$\\
 

   \begin{prop}\label{15-5}
Let $E$ be a geometric (a,b)-module with Bernstein polynomial $B_E$. Assume that $B_E$ has the root $-\alpha$ with multiplicity $p \geq 1$. Then there exists at least  $p$ distinct values of the integer  $j $ such that $-\alpha$ is a root of $B^j_E$.
\end{prop}

\parag{Proof} It is enough to proof the result when $E$ has a simple pole, by definition of the standard and  higher Bernstein polynomials.\\
We shall make an induction on the rank $r \geq 1$  of $E$. Since the case $r = 1$ is trivial, assume that the result is proved (for any $p \geq 1$) for each integer  $r \leq r_0$ and assume that the rank of $E$ is $r_0 + 1$. Then $G = E/S_1(E)$ has rank at most $r_0$ and we may apply the induction hypothesis to $G$ if $-\alpha$ is a root of order $p \geq 1$ of $B_G$.\\
Consider the exact sequence of geometric (a,b)-modules
$$ 0 \to S_1(E) \to E \to G \to 0 .$$
First assume that $B_{S_1(E)}(-\alpha) \not= 0$. Then the Bernstein polynomial of $G$ is divisible by $(x+\alpha)^p$ since we have an exact sequence of monodromic vector spaces
$$ 0 \to S_1(E)/bS_1(E) \to E/bE \to G/bG \to 0$$
corresponding to the previous exact sequence ($S_1(E)$ is normal in $E$). \\
In this case we conclude by the induction hypothesis applied to $G$ which has rank at most $r_0$ (since $E \not= \{0\}$ the rank of $S_1(E)$ is al least $1$). The conclusion follows from the remark following  point 6 (on $B_E$) recalled at the beginning of section 2.1.\\ 

Assume now that $B_{S_1(E)}(-\alpha) = 0$. Then, since $S_1(E)$ is semi-simple, $-\alpha$ is a simple root of $B_{S_1(E)}$ and the exact sequence above implies that 
$-\alpha$ is a root of order at least $p-1$ of $B_G$. Then the induction hypothesis and  Lemma \ref{21/1} give the existence of $(p-1)$ distinct  values of $j \geq 2$ such that  $B^j_E(-\alpha) = 0$. Since $B_{S_1(E)} = B^1_E$ we have found $p$ values of $j $ with $B^j_E(-\alpha) = 0$.$\hfill\blacksquare$\\

\parag{Warning}It happens that if $-\alpha$ is a root of multiplicity $p \geq 1$ of $B_E$ the $p$ values of $j$ for which $B^j_E(-\alpha) = 0$ do not contain $j = 1$ as is it shown by the following example.

\parag{Example} Define now the fresco  $E := B[a]\varphi_p \subset \Xi^{(N)}_\alpha$  with 
$$\varphi _p:= s^{\alpha+m-1}\frac{(Log s)^p}{p!} + s^{\alpha-1}$$
 where $m \geq1, p \geq 2$ and $\alpha \in ]0, 1[ \cap \mathbb{Q}$. Then the semi-simple part of $E^\sharp$ is $E_\alpha = Bs^{\alpha-1}$  and the Bernstein polynomial of $E$ is equal to 
  $(x+\alpha + m)^p(x + \alpha)$ since we have  $(b^{-1}a - (\alpha + m))[\varphi_p] = \varphi_{p-1} - (m+1)s^{\alpha -1}$. So $-\alpha - m$ is a root of $B^j_E$ for each $j \in [2, p+1]$ but not of $B^1_E(x) = x + \alpha$.  Note that $S_1(E) \simeq  E_{\alpha + p }$. $\hfill\square$\\

\begin{cor}\label{15-6}
Let $E$ be a geometric (a,b)-module with Bernstein polynomial $B_E$. Then $B_E$ divides  $ \prod_{j=1}^{d(E)} B^j_E(x) $ the product of the higher Bernstein polynomials $B^j_E$ of $E$.
\end{cor}

\parag{Proof} It is enough to consider the case of a simple pole $E$. We already know that each root of a polynomial  $B^j_E$ is a root of $B_E$ thanks to Proposition \ref{2} and the fact  that a root of multiplicity $p$ in $B_E$ is a root of at least $p$ polynomials $B^j_E$, thanks to Proposition \ref{15-5}. So $B_E$ divides the product $ \prod_{j=1}^{d(E)} B^j_E$.$\hfill\blacksquare$

\subsection{Complements}

We give now some more properties of the higher Bernstein polynomials of a geometric (a,b)-module which are useful in the sequel.

\begin{lemma}\label{30/1-b}
Let $G \subset H$ be two geometric  (a,b)-modules such that $H\big/G$ is a finite dimensional complex vector space. If $-\beta$ is a root of $B_G$  then there exists a root of $B_{E}$ in $-\beta + \mathbb{N}$.
\end{lemma}

\parag{Proof} We may assume that $G$ and $H$ has simple poles.\\
We shall prove the lemma by induction on the rank of $H$ (equal to the rank of $G$). In rank $1$ we have $G = E_\beta$ and $H = E_{\beta - \nu}$ with $\nu \in \mathbb{N}$, so the result is clear.\\
Assume that the result is proved for $rk(H) = k$ and consider  $G \subset H$ a sub-module of finite co-dimension in $H$  with the rank  of $H$ equal to $ k+1$. If $B_H$ has no root in $-\beta + \mathbb{N}$ then there exists a normal rank $1$ sub-module $E_\gamma$ in $S_1(H)$ and since the roots of $B^1_H$ are roots of $B_H$, we have $\gamma \not\in \beta - \mathbb{N}$. The sub-module $E_\gamma \cap G$ is isomorphic to $E_{\gamma + m}$ for some integer $m \geq 0$, and it is normal in $G$. So the Bernstein polynomial of  $G\big/E_{\gamma+m}$ has still the root $-\beta$ and since $G\big/E_{\gamma+m}$ is a sub-module with finite co-dimension in $H\big/E_\gamma$ the induction hypothesis gives the existence of a root in $-\beta + \mathbb{N}$ for the Bernstein polynomial of $H\big/E_\gamma$. Since this polynomial  divides $B_H$ we obtain a contradiction.$\hfill \blacksquare$\\

\begin{cor}\label{31/1}
Let $G \subset H$ be two geometric  (a,b)-modules such that $H\big/G$ is a finite dimensional complex vector space. If $-\beta$ is a root of $B^p_G$  then there exists $q \geq p$ such that  $B^q_H$ has a root  in $-\beta + \mathbb{N}$.
\end{cor}

\parag{Proof} We may assume that $G$ and $H$ has simple poles.\\
Let $\tilde{G} := G\big/S_{p-1}(G)$ and $\tilde{H}:= H\big/S_{p-1}(H)$. Then we have $B^1_{\tilde{G}} = B^p_G$ and $B^{p+h-1}_{\tilde{H}} = B^h_{\tilde{H}}$.
Since $B^p_G = B^1_{\tilde{G}}$ divides $B_{\tilde{G}}$ we may apply Lemma \ref{30/1-b} and we find that $B_{\tilde{H}}$ has a root in $-\beta + \mathbb{N}$ and so there exists $h \geq 1$ such that $B^h_{\tilde{H}} = B^{p+h-1}_H$ has such a root.$\hfill\blacksquare$\\

\begin{lemma}\label{30/1-a} Let $G \subset E$ be two geometric (a,b)-modules. Assume that $B_E$ has no root in $-\beta+\mathbb{N}^*$ and that $-\beta$ is 
a root of  $B^p_{G}$. Then $-\beta$ is a root of  a polynomial  $B^q_{E}$ for some $q \geq p$.
\end{lemma}

\parag{Proof} We may assume that $G$ and $E$ has a simple pole. Since $G$ has finite co-dimension in $N_E(G)$, using Lemma \ref{31/1} we may assume, replacing $G$ by $N_E(G)$,  that $G$ is normal in $E$ with the same hypothesis (but may be replacing $-\beta$ by $-\beta + m$ with $m$ a non negative integer).
 Since $B^p_{G}$ divides $B_{G/S_{p-1}(G)}$ and since we have  the exact sequence
$$ 0 \to G\big/S_{p-1}(G) \to E\big/ S_{p-1}(G) \to  E\big/\big(G + S_{p-1}(G)\big) \to 0 $$
 because we assume that $G$ is normal in $E$ (so $S_{p-1}(G)$ is normal in $E$), $B^p_G$ divides the Bernstein polynomial of  $E\big/S_{p-1}(G)$.\\
But this polynomial divides $B_{E\big/S_{p-1}(E)}$ which divides the product of the $B^j_{E\big/S_{p-1}(E)} $. 
Since we have $B^j_{E\big/S_{p-1}(E)} \simeq B^{j+p-1}_E$ we find that $-\beta$ must be a root of some $B^q_E$ for at least one integer  $q \geq p$.$\hfill \blacksquare$\\

The following simple lemma is proved in \cite{[part I]} Lemma 6.3.6

\begin{lemma}\label{30/1-d}
For any geometric (a,b)-module $E$ and any $\alpha \in \mathbb{Q}\big/\mathbb{Z}$ we have 
$$  \qquad \qquad \qquad  (E^{[\alpha]})^\sharp = (E^\sharp)^{[\alpha]} \qquad \qquad \qquad \qquad \qquad \qquad \qquad \qquad  \blacksquare$$
\end{lemma}

\smallskip

\begin{lemma}\label{30/1-c} Let $E$ be a  geometric (a,b)-module and assume that $B_E$ has a root $-\beta$ of multiplicity $p$. Then $-\beta$ is also a root of multiplicity $p$ in $B_{E^{[\alpha]}}$, where $[\alpha]$ is the class of $\beta$ in $\mathbb{Q}\big/\mathbb{Z}$.
\end{lemma}

\parag{Proof} We may assume that $E$ has a simple pole thanks to the previous lemma since for any geometric (a,b)-module $E$ we have $d(E) = d(E^\sharp)$.
Now the exact sequence 
$$ 0 \to E_{[\not= \alpha]} \to E \to E^{[\alpha]} \to 0 $$
splits and $B_E$ is the product of the Bernstein polynomials of $ E_{[\not= \alpha]}$ and of $E^{[\alpha]}$ thanks to the corresponding exact sequence of vector spaces
$$0 \to E_{[\not= \alpha]}\big/bE_{[\not= \alpha]} \to E\big/bE \to E^{[\alpha]}\big/bE^{[\alpha]} \to 0 $$
compatible with the respective actions of $-b^{-1}a$ and the fact that there is no eigenvalue of the action of $-b^{-1}a$ on $E_{[\not= \alpha]}\big/bE_{[\not= \alpha]}$ in $-\alpha + \mathbb{Z}$ so are disjoint from the eigenvalues of the action of $-b^{-1}a$ on $E^{[\alpha]}\big/bE^{[\alpha]}$.\\
Then the existence of the root  $-\beta \in -\alpha + \mathbb{Z}$ of multiplicity $p$ in $E/bE$ implies  the same property for $E^{[\alpha]}\big/bE^{[\alpha]}$.$\hfill \blacksquare$\\

\begin{lemma}\label{1/2-a}
Let $E$ be a simple pole geometric (a,b)-module. Then for each $[\alpha]$ in  $\mathbb{Q}\big/\mathbb{Z}$ and each integer  $j \in [1, d(E)]$ we have 
$$ B^j_E = B^j_{E_{[\not=\alpha]}}B^j_{E^{[\alpha]}}\quad {\rm and \ also } \quad  B_E = B_{E_{[\not=\alpha]}}B_{E^{[\alpha]}}.$$
\end{lemma}

\parag{Proof} This is obvious since we know that $E = E_{[\not=\alpha]} \oplus E^{[\alpha]}$ as a $B[a]$-module. Then for each $j$ we have
$S_j(E) =   S_j(E_{[\not=\alpha]}) \oplus S_j(E^{[\alpha]}).$\\
The conclusion follows because the eigenvalues of the  action of  $-b^{-1}a$ on the finite dimensional vector spaces $G\big/bG$ where $G$ is respectively 
$ S_j(E_{[\not=\alpha]}) \big/ S_{j-1}(E_{[\not=\alpha]}) $ and $S_j(E^{[\alpha]})\big/S_{j-1}(E^{[\alpha]})$ are mutually disjoint.$\hfill\blacksquare$\\

\begin{cor}\label{1/2-b}
Let $E$ be a geometric (a,b)-module and let $[\alpha]$ be in $\mathbb{Q}\big/\mathbb{Z}$. Then the polynomial $B_E$ and  for each integer $j \in [1, d(E)]$ the polynomial 
$B^j_{E^{[\alpha]}}$ are respectively obtained by  deleting in $B_E$ and  in $B^j_E$ all roots which do not belong to $[\alpha]$.
\end{cor}

\parag{Proof} Thanks to Lemma \ref{30/1-d} and Lemma \ref{1/2-a} we may assume that $E$ has a simple pole. 
Then Lemma \ref{1/2-a} gives the result.$\hfill \blacksquare$\\

 \section{The case of a fresco}
 
 \subsection{Some known facts}
 
First let me remind some previous results and fix some notations

\begin{enumerate}
\item  A {\bf fresco} is a geometric (a,b)-module $F$ such that $F\big/(aF + bF)$ is a one dimensional complex vector space. This is equivalent of the fact that $F$ is a sub-module of some $S\Xi^{(N)}_\mathscr{A}\otimes V$ generated by one element over $B[a]$.
\item It is proved in \cite{[B.09]}\footnote{The formal results of \cite{[B.09]} extend  to the convergent frescos. The division by $(a - \lambda b)$  in $B[a]$ is an easy exercise. } that any fresco $F$ is isomorphic to a quotient $B[a]\big/ B[a]P$ where $P$ is an element in $B[a]$ of the type
$$ P := (a -\lambda_1b)S_1^{-1}(a - \lambda_2b)S_2^{-1} \cdots (a - \lambda_kb) $$
where $k$ is the rank of $F$ as a $B$-module and where $-\lambda_j + k -j$ are negative rational numbers which are exactly the roots (counting multiplicities) of $B_F$.\\
Note that for any fresco $F$, $B_F$ is the {\bf characteristic polynomial} of the action of $-b^{-1}a$ on the vector space $F^\sharp\big/bF^\sharp$ which is equal to the minimal polynomial in this case.
\item If we have an exact sequence of geometric (a,b)-modules
 $$0 \to E_1 \to E \to E_2 \to 0$$
 and if $E$ is a fresco, then $E_1$ and $E_2$ are frescos and we have the equality
 $$ B_E(x) = B_{E_1}(x + r)B_{E_2}(x) $$
 where $r$ is the $B$-rank of $E_2$.
 \end{enumerate}
 
 \subsection{The $B^j_F$ when $F$ is a fresco}
 
 Our aim is now to prove that the definition given above for the higher  Bernstein polynomials of a geometric (a,b)-module is compatible with the definition of the higher  Bernstein polynomials of a fresco which is given in a previous version of this paper (see \cite{[B.23]}). For instance, we will verify that in the case of a fresco  $F$ the product of the $B^j_F$ is equal to $B_F$.\\
 The definition of the higher Bernstein polynomial of a fresco $F$ given in \cite{[B.23]} is the following.
 
 \begin{defn}\label{22/1-1}
 Let $F$ be a fresco and  let $j$ be an integer in $[1, d(F)]$. Then $B^{(j)}_F$, the $j$-th Bernstein polynomial of $F$ is defined by the formula
 \begin{equation*}
 B^{(j)}_F(x) = B_{S_j(F)\big/S_{j-1}(F)}(x + r_j)  \tag{$@$}
 \end{equation*}
 where $r_j$ is the rank of $F\big/S_{j+1}(F)$.
 \end{defn} 
 
 Our goal is to prove the following compatibility of this definition with the definition \ref{0} given in section 2.
 
 \begin{thm}\label{22/1-2}
 For each fresco $F$ and each integer $j$ we have $B^j_F = B^{(j)}_F$.
 \end{thm}
 
 This result will be the content of Corollary \ref{principal} below.

Recall that for any (a,b)-module $E$ with rank $r$ a {\bf Jordan-H\"{o}lder sequence}\footnote{Corollary 3.2.5 in \cite{[part I]} implies the existence of a normal sub-module of rank $1$. An induction is then enough to obtain the existence of a J-H. sequence for any regular (convergent)  (a,b)-module.} for $E$ is an increasing sequence of normal sub-modules
$$\{0\} = G_0 \subset G_1 \subset \dots \subset G_r = E$$
such that $G_{j+1}\big/G_j$ has rank $1$ for each $j \in [0, k-1]$. Any geometric (a,b)-module\footnote{The regularity of $E$ is enough in fact, as explained in the previous footnote.}
admits a Jordan-H\"{o}lder sequence.

\begin{prop}\label{16/12}
Let $F$ be a fresco with rank $r$ and let 
$$\{0\} = F_0 \subset F_1 \subset \dots \subset F_r = F$$
 be a Jordan-Hold\"er sequence for $F$. Then the normalization of $F_j^\sharp$ in $F^\sharp$ is equal to $b^{j-r}F_j^\sharp$. Moreover we have
$$ F^\sharp = \sum_{j=1}^r b^{j-r}F_j .$$
\end{prop}

Note that this equality shows that $b^{j-r}F_j^\sharp$ is a sub-module of $F^\sharp$ where  $r-j$ is the co-rank of $F_j$ in $F$.

\parag{Proof} Since our assertions are obvious for $r = 1$, assume that they are already proved for $r \geq 1$ and we shall consider the case where $F$ is of rank $r+1$.\\
So we assume that $\{0\} = F_0 \subset F_1 \subset \dots \subset F_r \subset F_{r+1} = F$ and that  we have $F_r^\sharp = \sum_{j=1}^r b^{j-r}F_j $. Then, assuming that $F/F_r \simeq E_\alpha$, let $x$ be a generator of $F$ (as a $B[a]$-module) whose image in $E_\alpha$ is the standard generator of $E_\alpha$. 

\parag{Claim 1} Then $y := (a -\alpha b)x $ is a generator of $F_r$.

\parag{Proof of the claim 1} As $x$ is a generator of $F$ it is enough to show that if $P \in B[a]$ satisfies $Px \in F_r$ then there exists $Q \in B[a]$ such that $Px =Qy$.\\
Make the division of $P$ by $(a - \alpha b)$ in $B[a]$. This gives $P = Q(a - \alpha b) + R$ with $R \in B$. Since $Px \in F_r$ is equivalent to $Pe_\alpha = 0$ in 
$E_\alpha = F\big/F_r$, we find that $Re_\alpha = 0$ in $E_\alpha$ and so $R = 0$. So we obtain $Px = Qy$ proving the claim.\\

This implies that $b^{-1}y$ is in $F^\sharp$ and so that $b^{-1}(F_r^\sharp) = \big(b^{-1}F_r\big)^\sharp$ is contained  in $F^\sharp$.\\
 Indeed we have in any simple pole (a,b)-module  and  for each integer $p \geq 0$ the equality $b(b^{-1}a)^p = (b^{-1}a - 1)^pb$.\\
Remark now that the surjective $B[a]$-linear map $\pi_0: F \to E_\alpha$ extends to a surjective $B[a]$-linear  map $\pi : F^\sharp \to E_\alpha$ because $E_\alpha$ has a simple pole. Since $F_r$ is the kernel of $\pi_0$ the kernel of $\pi$ is the normalization of $(F_r)^\sharp$ in $F^\sharp$ (the inclusion between the kernel  of $\pi$ and the normalisation of $(F_r)^\sharp$ is clear and both are normal sub-modules with the same rank).

\parag{Claim 2 } We have $F^\sharp = Bx + b^{-1}(F_r)^\sharp$ and  the kernel of $\pi$ is equal to $b^{-1}(F_r)^\sharp$. So the normalisation of $(F_r)^\sharp$ is $b^{-1}(F_r)^\sharp$.\\

Note that $b^{-1}(F_r)^\sharp$ is an (a,b)-sub-module of $F^\sharp$ but that  $Bx$ is not stable by $a$ (at least when $r \geq1$);  so $F^\sharp$ is the direct sum of $Bx$ and $b^{-1}(F_r)^\sharp$ as a $B$-module but not as a $B[a]$-module.

\parag{Proof of the claim 2} We know (see Claim 1)  that $y$ is a generator of $F_r$ and that $b^{-1}y$ is in $F^\sharp$, so the inclusion $Bx + b^{-1}(F_r)^\sharp \subset F^\sharp$ is clear. \\
Conversely, if $z$ is in $F^\sharp$ there exists  $ z_1, \dots, z_p  \in F $  such that $z = \sum_{j=1}^p (b^{-1}a)^jz_j$. Now, for each $j \in [1, p]$ write $z_j = S_j(b)x + t_j$ where  $S_j \in B$ and $t_j \in F_r$. \\
Then we obtain the opposite inclusion because   $(b^{-1}a)^jBx \subset Bx + b^{-1}(F_r)^\sharp$ is consequence of the formula, for $S \in B$:
$$ (b^{-1}a)S(b)x = b^{-1}\big(\alpha S(b)bx + S(b)y + b^2S'(b)x) \in Bx + b^{-1}By \subset Bx + b^{-1}(F_r)$$
and we have  $(b^{-1}a)^j t_j \in (F_r)^\sharp \subset b^{-1}(F_r)^\sharp$. This proves our first assertion.\\
 
Since the restriction of $\pi$ to $Bx$ is bijective, we conclude that $Ker(\pi) \subset b^{-1}(F_r)^\sharp$. But there exists an integer $n \geq 1$ such that $b^n(F_r)^\sharp$ is contained in $F_r$ and so $b^{n-1}(F_r)^\sharp$ is in  $ b^{-1}F_r$. This implies that $ b^{-1}(F_r)^\sharp  = (b^{-1}F_r)^\sharp \subset Ker(\pi)$ by $B$-linearity and normality of $Ker (\pi)$. So   Claim 2 is proved.\\
Now our induction hypothesis gives, since $Bx \subset F_{r+1}$ that $F^\sharp = \sum_{j=1}^{r+1} b^{j-r-1}F_j $.$\hfill\blacksquare$\\

We have proved in fact the more precise result for a fresco $F$ with rank $r$ and generator $x$:

\begin{cor} For a fresco $F$ with rank $r$ and generator $x$ there is   a direct sum decomposition of $F^\sharp$  as a $B$-module
$$ F^\sharp = \oplus_{j=1}^{r}  Bb^{j-r}x_j $$
where $x_j$ for each $ j \in [1, k]$ is the  generator of the  $F_j$ obtained  as follow:\\
Let $P$ be  a generator   in $B[a]$ of the left annihilator of $x$  in $F$  which may be written as (see point 2 above)
$$ P := (a - \alpha_1 b)S_1(a - \alpha_2 b)S_2 \dots S_{r-1}(a - \alpha_r b)S_r $$
where $S_1, \dots, S_r$ are invertible elements in $B$ and where the $x_j$ are defined by the formula
 $$x_j := (a - \alpha_{j+1} b)S_{j+1} \dots (a - \alpha_r b)S_r x \quad j \in [1, r-1] \quad  {\rm and} \quad x_r = x. \qquad \qquad \blacksquare $$.
 \end{cor}

\begin{cor}\label{18/12 0}
For a fresco $F$ and any normal sub-module $G$ of co-rank $g$ in $F$ the normalisation of $G^\sharp$ in $F^\sharp$ is equal to $b^{-g}G^\sharp$.
\end{cor}

\parag{Proof} With the notation of the previous proposition, we have shown that $b^{-1}(F_r)^\sharp$ is normal in $F^\sharp = (F_{r+1})^\sharp$. This implies, since each $F_h$ is normal and is  a fresco, that for each $h$ the sub-module  $b^{h-r-1}(F_h)^\sharp$ is normal in $F^\sharp = (F_{r+1})^\sharp$. Since it contains $(F_h)^\sharp$ and has the same rank, it is the normalization of $(F_h)^\sharp$ in $F^\sharp$.\\
Now the conclusion follows because we may assume that $G = F_h$ for some J-H. sequence of $F$: choose a J-H. sequence for $F/G$ , lift it in $F$ and complete with a J-H. sequence for $G$.$\hfill \blacksquare$\\

As a consequence we obtain that for any fresco $F$ and any $j\in [1, d(F)]$ we have
$$ S_j(F^\sharp) = b^{-r_j}S_j(F)^\sharp $$
where $r_j$ is the co-rank of $S_j(F)$ in $F$. So the    Bernstein polynomial  of $S_j(F^\sharp)$ is deduced from the Bernstein polynomial of $S_j(F)$ by the formula
$$ B_{S_j(F^\sharp)}(x) = B_{S_j(F)}(x + r_j) .$$

\begin{cor}\label{principal}
For each fresco $F$ and each integer  $j \in [1, d(F)]$, the equality
$$ B_{S_j(F^\sharp)/S_{j-1}(F^\sharp)}(x) =  B_{S_j(F)/S_{j-1}(F)}(x + r_j)$$
holds true. 
\end{cor}

\parag{Proof} We shall prove that there is an isomorphism of (a,b)-modules between $\big( S_j(F)/S_{j-1}(F)\big)^\sharp$ and $b^{r_j}\big(S_j(F^\sharp)/S_{j-1}(F^\sharp)\big)$. We have
\begin{equation*}
 \big(S_j(F)/S_{j-1}(F)\big)^\sharp = S_j(F)^\sharp \big/N_{S_j(F)^\sharp}(S_{j-1}(F)^\sharp) 
 \end{equation*}
  thanks to  Lemma  \ref{21/1}. We shall use now the following two easy properties for $X \subset Y$ an inclusion of regular  (a,b)-modules (see \cite{[part I]}):
  \begin{enumerate}
  \item For any non negative integer  $s $, $N_{Y}(b^sX) = N_Y(X)$.
  \item For any positive  integer  $r $, we have $b^{-r}X \subset b^{-r}Y$ and $N_{b^{-r}Y}(b^{-r}X) = b^{-r}\big(N_Y(X)\big)$.
  \item For any positive  integer $j$ and any $n \in \mathbb{Z}$ we have $S_j(b^nX) = b^nS_j(X)$.
  \end{enumerate}
  The proof is left to the reader.\\
  Denote by $r_j$ the co-rank of $S_j(F)$ in $F$ and by $s_j$ the co-rank of $S_{j-1}(F)$ in $S_j(F)$. Then we have  $S_j(F^\sharp) = b^{-r_j}S_j(F)^\sharp$ and
  $S_{j-1}(F^\sharp) = b^{-r_j-s_j}S_{j-1}(F)^\sharp$ thanks to Corollary \ref{18/12 0}. Using the three  properties  above we see that
  $$N_{S_j(F)^\sharp}(S_{j-1}(F)^\sharp) = b^{r_j}\big(N_{S_j(F^\sharp)}(b^{s_j}S_{j-1}(F^\sharp)) = b^{r_j}N_{S_j(F^\sharp)}(S_{j-1}(F^\sharp))  = b^{r_j}S_{j-1}(F^\sharp) $$
  since $S_{j-1}(F^\sharp)$ is normal in $S_j(F^\sharp)$. Since $S_j(F)^\sharp = b^{r_j}S_j(F^\sharp)$ the proof is complete using the easy fact that for $X$ a normal sub-module in an (a,b)-module $Y$ we have for each integer $r \in \mathbb{Z}$ a natural isomorphism 
  $ b^r(Y/X) \simeq b^rY/b^rX.\hfill \blacksquare$\\
  
  As a consequence, the definition of the polynomials $B^j_F$ given in \cite{[B.23]} for a fresco $F$ coincide with the general definition of the higher Bernstein polynomials of a geometric (a,b)-module given in  the previous section.

 \begin{cor}\label{20/1}
 For a fresco $F$ we have the equality
 $$B_F = \prod_{j=1}^{d(F)} B^j_F .$$
 \end{cor}
 
 \parag{Proof} We know that $B_F$ and the product $\prod_{j=1}^{d(F)} B^j_F$ have the same degree equal to the rank of $F$. Since $B_F$ divides $\prod_{j=1}^{d(F)} B^j_F $ and  they are both monic, the equality follows.$\hfill \blacksquare$\\
 
  \subsection{Complements}
 
 Since in our Theorem \ref{1} we consider a fresco $\F$ with $d(\F^{[\alpha]}) = d$ and a root $-\alpha-m$ of $B^d_\F$ we want to show that we may reach this situation if we begin with the hypothesis that $B^d_{H^{n+1}_0}$  has a root $-\alpha -m'$. The proposition \ref{5/1/24} below shows that in this situation we may find such a fresco $\F$ inside $H^{n+1}_0$ with $B^d_\F$ having a root $-\alpha -m$ with $m \geq m'$. In the case where $d = d(\h^{[\alpha]})$ we have moreover $m = m'$. \\
  
 \begin{prop}\label{5/1/24}
Let $\E$ be an $[\alpha]$-primitive geometric (a,b)-module and let $j$ be an integer in $[1, d(\E)]$ where $d(\E)$ is the  nilpotent order   of $\E$ . Assume that $ -\alpha - m$ is a root of the Bernstein polynomial $B^j_{\E}$. \\
Then there exist an integer $m' \geq m$ and  a fresco $\F \subset \E$  with that $d(\F) \geq j$  such that $-(\alpha + m')$ is a root of some  $h$-th Bernstein polynomial  $B^h_{\F}$ of $\F$, for some $h \geq j$.\\
Moreover, if $j = d(\E)$ then we may take $m' = m$ with  $h = d(\E)$.
\end{prop}

Note that, thanks to Corollary \ref{23/12}, we have also  a root $-(\alpha + m')$ for $B^j_{\F}$  but may be with $m' \in \mathbb{Z}$.

\parag{Proof} Our hypothesis implies that there exists a $B[a]$-linear surjective map, where $\beta := \alpha + m$: 
$$ \pi :  S_j(\E^\sharp)/S_{j-1}(\E^\sharp) \to E_\beta .$$
Let $n \in \mathbb{N}$ such that $b^nS_j(\E^\sharp) \subset S_j(\E)^\sharp$. Remark that if
 we assume $j = d(\E)$ we have $S_j(\E^\sharp) = \E^\sharp$ we may take $n = 0$.\\
 We also denote by $\pi$ the map  $S_j(\E^\sharp) \to \E_\beta$ obtained by composition with the quotient map to $S_j(\E^\sharp)/S_{j-1}(\E^\sharp)$.\\
Take $x \in S_j(\E^\sharp)$  such that $\pi(x) = e_\beta$ where $e_\beta$ is the standard generator of $\E_\beta$ (so $\E_\beta = B.e_\beta$ and $a e_\beta = \beta b e_\beta$). Since there exists some non negative integer $n$ such that $b^nx$ is in $S_j(\E)^\sharp$ we may find $y_0, \dots, y_p$ in $S_j(\E)$ such that $ b^n x  = \sum_{p=0}^{N} (b^{-1}a)^p y_p $.  Then there exists at least one $p_0 \in [0, N]$ such that $\pi(y_p)$ is not contained in $b^{n+1}\E_\beta$. Since $\pi$ vanishes on $S_{j-1}(\E^\sharp)$ we have $y_{p_0} \not\in S_{j-1}(\E)$. Then applying a suitable invertible element in $B$ to $y_{p_0}$ we obtain an element $z \in S_j(\E) \setminus S_{j-1}(\E)$ such that the  fresco $\F := B[a]z$  satisfies $d(\F) \geq j$ and $\F\big/S_{j-1}(\F)$ has a surjective map to $\E_{\alpha +m'}$ with $m' \leq m + n$. So $-(\alpha+m')$ is a root of  the Bernstein polynomial of $\F\big/S_{j-1}(\F)$ and so of at least one polynomial $B^h_{\F\big/S_{j-1}(\F)} = B^{j+h-1}_\F$ for some $h \geq 1$. The proof is complete.$\hfill \blacksquare$\\

The following example shows that, in general, for $1 \leq j < d(\E)$ we cannot avoid the integer shift for the root of the Bernstein polynomial of the fresco constructed in the previous proposition.

\parag{Example}  Fix $\beta \not= \gamma$ in $]0, 1[ \cap \mathbb{Q}$. Let $\E$ be the geometric (a,b)-module generated by $\xi := s^{\beta-1}Log\, s \otimes v_1$ and $\eta := s^{\gamma-1} \otimes v_2$ inside $\Xi_{\beta, \gamma}^{(1)} \otimes V$ where $(v_1, v_2)$ is a basis of $V$. Then we have $S_2(\E) = \E$, so $d(\E) = 2$,
$$ \E^\sharp = \E + \C s^{\beta-1}\otimes v_1, \quad  S_1(\E) = S_1(\E)^\sharp = (\E_{\beta+1}\otimes v_1) \oplus (\E_\gamma \otimes v_2) $$
and 
$$ S_1(\E^\sharp) = (\E_\beta \otimes v_1) \oplus (\E_\gamma \otimes v_2).$$
So $-\beta$ is a root of $B^1_\E$.\\
The only surjective map $\pi : S_1(\E^\sharp) \to \E_\beta$ is  given (up to a multiplicative constant) by the projection of $ S_1(\E^\sharp)$ onto $\E_\beta \otimes v_1$ and the only choice for an element  $x$ such that  $\pi(x) = e_\beta$  is $ x = s^{\beta-1}\otimes v_1$ up to  a multiplicative constant and an element in  $ \E_\gamma \otimes v_2$ . Since we have $x = -(b^{-1}a - \beta)\xi $ modulo $\E_\gamma\otimes v_2$ the map $\pi$ is not defined on $\xi$ but only on $bx$ which is in $S_1(\E)^\sharp  = S_1(\E)$. So the fresco $B[a]bx$ which is semi-simple has $-(\beta+1)$ (and may be \,  $-\gamma$) as root of its Bernstein polynomial. There is no fresco in $\E$ which is semi-simple and having $-\beta$ for root of its Bernstein polynomial.\\

\subsection{Jordan blocs}

Recall that a {\bf theme} is, by definition, a fresco which may be embeded in some $S\Xi_{\mathscr{A}}^{(N)}$ (so with $V = \mathbb{C}$) (see Definition 5.1.3 in \cite{[part I]}).

\begin{lemma}\label{Jordan bloc dans theme}
Let $\alpha \in ]0,1[ \cap \mathbb{Q}$ and let $\varphi$ be an element in $S\Xi_\alpha^{(N)}$ which is of degree $N$ in $Log\, s$. Then inside the rank $N+1$ theme $ T := B[a].\varphi \subset S\Xi_\alpha^{(N)}$, there exists an element 
 $$ \psi_N := s^{\alpha +m-1} \frac{(Log\, s)^N}{N!}$$
where $m$ is an integer.\\
For $\alpha = 1$, if $\varphi \in S\Xi_1^{(N)}$ is the class of a series which has degree $N+1$ in $Log\, s$,  then inside the rank $N+1$ them $ T := B[a]\varphi \subset S\Xi_1^{(N)}$  there exists an element which is the class of $s^m \frac{(Log\, s)^{N+1}}{(N+1)!}$, where $m$ is an integer.
\end{lemma}

\parag{Proof} Note first that $T := B[a] \varphi$ is a rank $(N+1)$ fresco  thanks to Lemma 5.2.4 in \cite{[part I]}.\\
  We shall prove the lemma for $\alpha \not= 1$ by induction on $N \geq 0$. We leave the case $\alpha = 1$ which is similar as an exercise for the reader.\\
   Since the case $N = 0$ is clear, assume that the lemma is proved for $N-1$ and let $T \subset S\Xi_\alpha^{(N)}$ a rank $N+1$ fresco. Then $S_1(T)$ is equal to $T \cap S\Xi_\alpha^{(0)}$ and we may embed the rank $N$ fresco $T/S_1(T)$ in
 $$S\Xi_\alpha^{(N-1)} \simeq S\Xi_\alpha^{(N)}/S\Xi_\alpha^{(0)}.$$ 
Thanks to our inductive hypothesis there exists an integer $m'$ such that 
$$ s^{\alpha+m'-1}(Log\, s)^N/N! \qquad  modulo \  S\Xi_\alpha^{(0)}$$
 is in $T/S_1(T)$ and, since $S_1(T) \subset S\Xi_\alpha^{(0)}$,  we may find an invertible element $S$ in $B$ such that $\varphi := s^{\alpha+m'-1}(Log\, s)^N/N! + S(b)s^{\alpha + M-1}$ is in $T$. Since $S_1(T)$ is isomorphic to $\mathcal{E}_{\alpha+q}$ for some positive integer $q$, for an integer $m''$  large enough, $s^{\alpha+m'+m''-1}(Log\, s)^N/N! $ will be in $T$, since $S(b)s^{\alpha+M+m''-1}$ will be in $S_1(T)$, concluding the proof.$\hfill \blacksquare$

\begin{cor}\label{Bloc de Jordan dans fresque}
Let $\mathcal{F}$ be a fresco and assume that the $p$-th Bernstein polynomial of $\mathcal{F}$ has a root in $-\alpha -\mathbb{N}$, where $\alpha$ is in $]0,1] \cap \mathbb{Q}$. Then there exists $w_1, \dots, w_p$ in $\mathcal{F}$ (in fact in $\mathcal{F}_{[\alpha]}$)  and an integer $m \in \mathbb{N}$ satisfying the relations:
\begin{equation*}
	aw_j = (\alpha +m)bw_j + bw_{j-1} \quad \forall j \in [1,p] \quad {\rm with \ the \ convention} \quad w_0 \equiv 0 \tag{$\star$}
\end{equation*}
and which are $B$-linearly independent in $\mathcal{F}$.
\end{cor}

\parag{Proof} We assume that $\alpha \not= 1$ leaving the case $\alpha = 1$ which is analogous  as an exercise for the reader.\\
Since $S_p(\mathcal{F}_{[\alpha]}) = S_p(\mathcal{F})_{[\alpha]}$, thanks to Lemma 4.2.5 in \cite{[part I]},
we may find an $[\alpha]$-primitive theme $T$ of rank $p$ in $\mathcal{F}_{[\alpha]}$ thanks to Proposition 6.3.3 in \cite{[part I]}. 
As we may assume that $T$ is embedded in $S\Xi_\alpha^{(p-1)}$ the previous lemma shows that there exists an integer $m_0$ such that $s^{\alpha+m_0-1}(Log\, s)^p/p!$ is an element in $T$.\\
Define $w_j = s^{\alpha+m_0-1}(Log\, s)^j/j!$ for $j \in [1,p]$. Then the relations $(\star)$ are satisfied and imply that $w_1, \dots, w_p$ are elements in $T^\sharp$.\\
To shows that $w_1, \dots, w_p$ are $B$-linearly independent, note  $J$ the $B$-sub-module generated by $w_1, \dots, w_p$. Then it has rank at most $p$. But the relation $(\star)$ shows that $J$ is an (a,b)-sub-module of $T^\sharp$ with a simple pole. Moreover as $w_1, \dots, w_p$ are clearly linearly independent over $A = \mathbb{C}\{s\}$, we have
$$ \dim_{\mathbb{C}} J/aJ = \dim_{\mathbb{C}} J/bJ \geq p $$
so $J$ has rank $p$ as a $B$-module.\\
Since there exists a non negative integer $q$ such that $b^qT^\sharp \subset T \subset \mathcal{F}_{[\alpha]}$ the conclusion follows considering $b^qw_1, \dots, b^q w_p$ and $m := m_0 + q$. $\hfill \blacksquare$\\

\parag{Remark} Let $J := \sum_{j=1}^p Bw_j$ the sub-$B$-module generated by $w_1, \dots, w_p$. Then $J$ is a sub-(a,b)-module of $\mathcal{F}$ which has a simple pole and is
 $[\alpha]$-primitive; it is equal to $\mathcal{E}(J_{\alpha+m, p})$ where $J_{\alpha+m, p}$ is the matrix of the standard Jordan bloc with rank $p$ and eigenvalue $\alpha+m$ (see  the end of Section 2.3 in \cite{[part I]}). The action of $b^{-1}a$ on $J/bJ$ is given by $J_{\alpha+m, p}$. So the Bernstein polynomial of $J$ is equal to $(x+ \alpha+m)^p$.\\

It is interesting to compare this result with Corollary 3.2.5 in \cite{[part I]}. Here we do not assume that the action of $b^{-1}a$ on $\mathcal{F}^\sharp/b\mathcal{F}^\sharp$ has a Jordan block of size $p$ for some $\lambda$ in $\alpha + \mathbb{N}$ but, in a way, that this happens for the eigenvalue $\exp(2i\pi\alpha)$ of   $\exp(2i\pi b^{-1}a)$ acting of $\mathcal{F}^\sharp/b\mathcal{F}^\sharp$. And this hypothesis  is precisely formulated by the existence of a root in $-\alpha - \mathbb{N}$ for the $p$-th Bernstein polynomial of the fresco $\mathcal{F}$.\\
Note that contrary to the result in Corollary 3.2.5 in  \cite{[part I]}  we have no control  here on the  integral shift between the root of the $p$-th Bernstein polynomial and the (multiple) root of the Bernstein polynomial  of the Jordan block obtained.\\

  \section{Existence of poles}
   
    \subsection{The complex of sheaves $(Ker\, df^\bullet, d^\bullet)$}
    
    \parag{The standard situation } 
    
     We consider now the following situation:\\
     Let  $f : U \to \mathbb{C}$ be a holomorphic function on an open polydisc $U $ with center $0$ in $\mathbb{C}^{n+1}$. We shall assume that $U$ is small enough in order that the inclusion $\{ df = 0 \} \subset \{f = 0\}$ holds in $U$.\\ 
     We denote $Y$ the hypersurface $\{ f = 0\}$ in $U$ and we assume that $Y$ is reduced. For each point $y \in Y$ we denote $f_y : X_y \to D_y$ a Milnor representative of the germ of $f$ at $y$. So $X_y$ is constructed by cutting a small ball, with center $y$ and with radius $\varepsilon > 0$ very small, with $f^{-1}(D_\delta)$ where $D_\delta$ is an open disc with center $0$ and radius $\delta \ll \varepsilon$. For $y = 0$ we simply write $f : X \to D$ such a Milnor representative of the germ of $f$ at the origin.\\
     Let $\pi: H \to D^*$ be the universal cover of the punctured disc $D^* := D \setminus \{0\}$ and choose a base point $\tilde{s}_0$ in $H$ over the base point $s_0$ in $D^*$. Fix a point $y \in Y$ and take for $D$ the disc of a  Milnor representative of $f_y$. Then we identify the Milnor fiber $F_y$ of $f$ at $y$ with $f^{-1}(s_0)$.\\
  For any $p$-cycle $\gamma$ in $H_p(F_y, \mathbb{C})$ let $(\gamma_{\tilde{s}})_{\tilde{s} \in H}$ be the horizontal family of $p$-cycles in the fibers of $f\times_{D^*}\pi$ over $H$ taking the value $\gamma$ at the point $\tilde{s}_0$. Then the regularity of the Gauss-Manin connection of $f$ at $y$  insures that for any $\omega \in \Omega_y^{p+1}$ which satisfies $d\omega = 0 $ and $df\wedge \omega = 0$ the (multi-valued) function $s \mapsto \int_{\gamma_s} \omega/df$ has a convergent asymptotic expansion when $s$ goes to $0$, which is in $\Xi_{\mathscr{A}}^{(p-1)}$ where $\exp(2i\pi\mathscr{A})$ contains the eigenvalues of the monodromy of $f$ at the point $y$.\\
     We define on $Y$ the following  sheaves ( see \cite{[B.I]} or \cite{[B.II]} ) for each integer $p \in [1,n]$:\\
     First let $Ker\, df^{p+1} \subset \Omega^{p+1}$ be the kernel of the  map $\wedge df : \Omega^{p+1} \to \Omega^{p+2}$ of coherent sheaves on $U$ and $Ker\, d^{p+1}$ be the kernel of the ($\mathbb{C}$-linear) de Rham  differential
         $$d^{p+1} : \Omega^{p+1} \to \Omega^{p+2}.$$
      Then for $p \in [1,n]$ define the sheaf $\mathcal{H}^{p+1}$  as the (topological) restriction on $Y$ of  the sheaf $Ker\, df^{p+1} \cap Ker\, d^{p+1}\big/ d(Ker\, df^{p})$.\\
      By convention we put $\mathcal{H}^{p+1} = 0$ for $p \not\in [1,n]$.\\
     Then we have a natural structure of $A$-modules on the sheaves $\mathcal{H}^{p+1}$ for each  $p$ induced by the natural action of $A$ on $\Omega_{\vert Y}$ given by $(g,\omega) \mapsto f^*(g) \omega$ where $g$ is in  $ A := \mathbb{C}\{s\}$ and $\omega$ is in  $ \Omega_y^{p+1}$, for each $y \in Y$.\\
     We have also an action of $\mathbb{C}[b]$ on $\mathcal{H}^{p+1}$ for each $p \in [1,n]$ which is defined as follows:
     \begin{itemize}
     \item For $\omega_y \in Ker \, d^{p+1} \cap Ker\, df^{p+1}$ write $\omega_y := du_y$ for some $u_y \in \Omega^p_y$ (holomorphic de Rham Lemma) and put $b[\omega_y] := [df\wedge u_y]$. \\
     Then clearly $d(df\wedge u_y) = 0$ and $df\wedge (df\wedge u_y ) = 0$.
     \item Of course, if we change the choice of $u_y$ ( for $p \in [1, n]$) in $u_y + dv_y, v_y \in \Omega_y^{p-1}$, the class of $b[\omega_y] \in \mathcal{H}^{p+1}$ is the same since  $df\wedge dv_y = - d(df\wedge v_y)$ is in    $d(Ker\, df ^p)$.  
       \end{itemize}
            
       The sheaf $\mathcal{H}^{p+1}$ modulo its $a$-torsion, noted $H^{p+1}$,   is the (a,b)-module version of the Gauss-Manin connection in degree $p$. As we assume $f$ reduced, the $0$-th cohomology of the Milnor fiber is $\mathbb{C}$ and the corresponding monodromy is trivial. 
     
     \begin{lemma}\label{(a,b) 1}
     The actions of $a$ and $b$ on $\mathcal{H}^{p+1}$ satisfy the commutation relation  \\ $ab - ba = b^2$.
     \end{lemma}
     
     \parag{Proof} For $\omega_y = du_y \in Ker\, df^{p+1} \cap Ker\, d^{p+1}$ we have $$b(a[\omega_y] + b[\omega_y]) = b[fdu_y + df\wedge du_y] = b[d(fu_y)] = [df\wedge fu_y] = ab[\omega_y]$$
     which gives the relation $b(a+b) = ab$ concluding the proof.$\hfill \blacksquare$\\
    
  Note that the action of $a$ is well defined on $Ker\, df^{p+1}$ but the action of $b$ is only well defined on the cohomology $\mathcal{H}_y^{p+1}$ for each $p \in [1,n]$ and each $y \in Y$.\\

\begin{thm}\label{embedding ex.}
We keep the notations introduced above and let $p$ be an integer in $[1,n]$. 
 Let $\omega \in \Omega^{p+1}_y$ be in $Ker\, df^{p+1} $ such that $d\omega = 0$. Then for each $\gamma \in H_p(F_y, \mathbb{C})$ define $\Phi(\gamma, \omega)$ as the element in $S\Xi_{\mathscr{A}}^{(p-1)}$ given by the singular part of the  asymptotic expansion
 \footnote{Since for the eigenvalue $1$ we consider only the {\em singular part} of the asymptotic expansion, so we replace  \ $\Xi_1^{(p-1)}$ by  \  
 $S\Xi_1^{(p-1)} := \Xi_1^{(p)}\big/ \Xi_1^{(0)}$ which is isomorphic to $\Xi_1^{(p-1)}$  ( see \cite{[B.84-b]}),  the
  $\Gamma$-factor that we introduce below shifts the order of poles at points in $-\mathbb{N}$ in the complex Mellin transform  $F^{\omega,\omega'}_h(\lambda)$  (see \cite{[B-M.89]}) of the associated  hermitian periods $\int_{f= s} \rho \omega\wedge \bar\omega'\big/df\wedge d \bar f$.}
  of the period-integral $\int_{\gamma_s} \omega/df$, where $\mathscr{A}$ is the image in $]0, 1]$ of the opposite of the roots of the reduced Bernstein polynomial of $f$ at the point $y \in Y$ and where $(\gamma_s)_{s \in H}$ is the horizontal family of $p$-cycles taking the value $\gamma$ at the base point $\tilde{s}_0$.  So we have:
$$ \Phi(\omega,\gamma) := \int_{\gamma_s} \omega/df \in S\Xi_{\mathscr{A}}^{(p-1)} .$$
Then using the fact that $H^p(F_y, \mathbb{C})$ is the dual of $H_p(F_y, \mathbb{C})$ and the linarity of $\Phi$ in $\gamma$,  $\Phi$ defines a map
$$ \Psi : \mathcal{H}^{p+1}_y \to S\Xi_\mathscr{A}^{(p)}\otimes_{\mathbb{C}} H^p(F_y, \mathbb{C}) ,\quad \Psi(\omega) := [ \gamma \mapsto \Phi(\omega,\gamma)] $$
which is $A$-linear and $b$-linear and whose kernel is equal to the $a$-torsion of $\mathcal{H}^{p+1}$.
\end{thm}

\parag{Proof} The $A$-linearity of $\Psi$ is obvious. The $b$-linearity is an easy consequence of the derivation formula
$$ \partial_s(\int_{\gamma_s} u) = \int_{\gamma_s} du/df $$
when $u$ is in $\Omega_y^p$ satisfies $df\wedge du = 0$.\\
Consider now $\omega \in Ker\, df^{p+1}$ such that $d\omega = 0$ and assume that $\omega$ is  in the Kernel of $\Psi$.  Then for each $\gamma$ the corresponding  period-integral vanishes because the asymptotic 
expansion is convergent (thanks to the regularity of the Gauss-Manin connection). So the class induced by $\omega/df$ in $H^p(F_y, \mathbb{C})$ vanishes which implies that the class defined by $\omega$
in the $f$-relative de Rham cohomology vanishes and so we may find a meromorphic form $ v \in \Omega_y^{p-1}[f^{-1}]$ such that $\omega = df\wedge dv$ (see \cite{[B.84-a]} and \cite{[B.84-b]} for $\alpha = 1$). This implies that $a^N[\omega] = 0 $ in $\mathcal{H}^{p+1}_y$.$\hfill \blacksquare$

\parag{Remark} The map $\Phi$ satisfies also the relation $\Phi(\omega, T\gamma) = \mathscr{T}(\Phi(\omega,\gamma))$ where $T$ is the monodromy acting on $H_p(F_y, \mathbb{C})$ and where $\mathscr{T}$ is the monodromy acting on $\Xi_{\mathscr{A}}^{(p-1)}$ \  via  \ $Log\, s \mapsto Log\, s +2i\pi$. So the image of $\Psi$ is contained in the sub-(a,b)-module of $S\Xi_{\mathscr{A}}^{(p-1)}\otimes_{\mathbb{C}} H^p(F_y, \mathbb{C})$ which is invariant by $\mathscr{T}\otimes T^*$ where $T^*$ is the action of the monodromy on $H^p(F_y,\mathbb{C})$.$\hfill \square$

\begin{cor}\label{fund. 1}
For $y \in Y$ let $H^{p+1}_y$ be the quotient of $\mathcal{H}^{p+1}_y$ by its $a$-torsion. Then $H^{p+1}_y$ is a geometric (convergent)  (a,b)-module.
\end{cor}

\parag{Proof} The point is that $H^{p+1}_y$ is a finite type $A$-module since $A$ is noetherian and  \ $S\Xi_{\mathscr{A}}^{(p)}$ is a finite type (free) $A$-module. Then $H^{p+1}_y$ is closed for the natural dual Fr\'echet topology induced by $ \Xi_\mathscr{A}^p\otimes H^p(F_y, \mathbb{C}) $. As it is $b$-stable it is also stable by the action of $B[a]$.  So it is a geometric (a,b)-module.$\hfill \blacksquare$\\

Note that it is not obvious  to show directly that $B$ acts on $H^{p+1}_y$  contrary to the formal case.

\begin{defn}\label{fresco of form} In the situation above, let $\omega$ be a germ at $y \in Y$ of the sheaf  $ Ker\, df^{p+1}_y$ which is $d$-closed. Then we define the {\bf fresco $\mathcal{F}_{f,\omega,y}$} associated to these data as the fresco $B[a] [\omega] \subset H^{p+1}_y$ which is generated in the geometric (a,b)-module  $H^{p+1}_y$ by the class of $\omega$.\\
\end{defn}

Note that for $p = n$ each germ $\omega$ at a point $y$ of $\Omega^{n+1}_y$ satisfies $df \wedge \omega = 0$ and $d\omega = 0$. \\
In the sequel we shall mainly use the case $p = n$ with  $y = 0$. So we simplify the notation to $\mathcal{F}_\omega$ when we consider the fresco $\mathcal{F}_{f,\omega, 0}$ in $H^{n+1}_0$.

 \subsection{The use of frescos} 
 
We begin by the  definition of the main hypothesis on the holomorphic germ $f$  which is assumed in the sequel.
 
 \begin{defn} \label{H} In the standard situation, fix a rational number $\alpha \in ]0, 1]$. We say that the germ $f$ has {\bf an isolated singularity for the eigenvalue $\exp(2i\pi\alpha)$ of its monodromy} when the local monodromy of $f$ at each point $y \not= 0$ in the reduced hypersurface $Y = f^{-1}(0)$,  acting on the {\em reduced} cohomology of the Milnor fiber at the point $y$ does not admit this eigenvalue. This hypothesis is denoted $H(\alpha, 1)$ in the sequel.
 \end{defn}
 Let us recall some known facts.
 \begin{enumerate}
 \item The hypothesis $H(\alpha, 1)$ is equivalent to the fact that, in  open neighborhood of the origin, the local {\em reduced}  $b$-function of $f$ at any point $x \not= 0$ has no root in $-\alpha - \mathbb{N}$.
 \item The hypothesis $H(\alpha, 1)$ is equivalent to the fact that, in an open neighborhood of the origin, the polar parts  of the meromorphic extension of the distributions 
 $$\frac{1}{\Gamma(\lambda)} \vert f\vert^{2\lambda} \bar f^{-h}, \quad \forall h \in \mathbb{Z},$$
  at points in $-\alpha - \mathbb{N}$ are distributions with support $\{0\}$.
 \item  The hypothesis $H(\alpha, 1)$ is equivalent to the fact that, for any test form $\varphi$ in $ \mathscr{C}^\infty_c(\mathbb{C}^{n+1})^{n+1, n+1}$ with compact support in $X \setminus \{0\}$  the meromorphic extension of the functions 
 $$ \frac{1}{\Gamma(\lambda)} \int_X \vert f\vert^{2\lambda} \bar f^{-h} \varphi $$
 has no pole in $-\alpha - \mathbb{N}$ for each $h \in \mathbb{Z}$.
 \item Since the monodromy of $f$ is  defined on $H^p(F_y, \mathbb{Z})$, for $\alpha \in ]0, 1[$ the hypothesis $H(\alpha, 1)$ is equivalent to the hypothesis $H(1-\alpha, 1)$. \\
 So the hypotheses of isolated singularity at the origin for the eigenvalues $\exp(\pm 2i±\pi\alpha)$ of the monodromy  for a holomorphic germ $f$ are equivalent.
 \end{enumerate}
 
 Assume that we are in the standard situation and that $f$ satisfies the hypothesis $H(\alpha, 1)$, that is to say that $f$ has an isolated singularity for the eigenvalue $\exp(2i\pi\alpha)$ of its monodromy.\\
Let $\omega, \omega' \in \Omega^{n+1}_0$ and let $\rho \in\mathscr{C}^\infty_c(\mathbb{C}^{n+1})$ such that $\rho \equiv 1$ near $0$ and its support is  small enough in order that $\rho \omega\wedge \bar \omega'$ is a well defined and $\mathscr{C}^\infty_c$ differential form of type $(n+1, n+1)$ on $\mathbb{C}^{n+1}$. Then for any $h \in \mathbb{Z}$ the holomorphic function, defined for $2\Re(\lambda) >  \sup \{0, h\}$ by the formula
\begin{equation*}
F^{\omega,\omega'}_h(\lambda) := \frac{1}{\Gamma(\lambda)} \int_X \vert f\vert^{2\lambda} \bar f^{-h} \rho \omega\wedge \bar\omega'  \tag{F}
\end{equation*}
has a meromorphic continuation to the all complex plane with poles contained  in $-\mathscr{A} - \mathbb{N}$ where $-\mathscr{A}$ is the finite subset of $\mathbb{Q}^-$  which is the set of the roots of the reduced Bernstein polynomial $\tilde{b}_{f, 0}$ of $f$ at the origin. \\
Moreover, thanks to our hypothesis $H(\alpha, 1)$  we have the following properties (see  \cite{[B.22]} for a proof) :
\begin{enumerate}
\item The polar parts of $F^{\omega, \omega'}_h(\lambda)$ on the points in $-\alpha - \mathbb{N}$ do not depend on the choice (with the conditions specified above) of the function $\rho$.
\item  The polar parts of $F^{\omega, \omega'}_h(\lambda)$ at points in $-\alpha - \mathbb{N}$ depend,  for given $\omega'$ and $h$, only on the image of $\omega$ in the formal (a,b)-module $\widehat{H}^{n+1}_0$, which is the formal completion of the geometric (a,b)-module $H^{n+1}_0$ defined in section 4.1.
\end{enumerate}

   The following result is proved in \cite{[B.22]} Proposition 3.1.1.
   
 \begin{prop}\label{tool 1} 
In the standard situation  assume  that the hypothesis $H(\alpha,1)$ is satisfied. Let $\omega$ and $\omega'$  be  holomorphic $(n+1)$-differential forms on $ X$  and let $\rho$ be a $\mathscr{C}^\infty$ function with compact support in  $ X$  which satisfies $\rho \equiv  1$ near $0$.  We have the following properties:
\begin{enumerate}[i)]
\item If there exists  $v \in \Omega^n(X)$ satisfying $df \wedge v \equiv 0$ and $dv = \omega$ on $X$, then $F^{\omega,\omega'}_{h}(\lambda)$ has no pole in $-\alpha - \mathbb{N}$ for any $h \in \mathbb{Z}$ and any $\omega' \in \Omega^{n+1}_0$.
\item $F^{a\omega, \omega'}_{h}(\lambda) - (\lambda+1)F^{\omega, \omega}_{h-1}(\lambda+1)$ has no pole in $-\alpha - \mathbb{N}$ for any $h \in \mathbb{Z}$ and any $\omega' \in \Omega^{n+1}_0$.
\item $F^{b\omega, \omega'}_h(\lambda) + F^{\omega, \omega'}_{h-1}(\lambda+1) $ has no pole in $-\alpha - \mathbb{N}$ for any $h \in \mathbb{Z}$ and any $\omega' \in \Omega^{n+1}_0$.
\item For any complex number $\mu$, $F^{(a - \mu b)\omega, \omega'}_h(\lambda) - (\lambda+\mu +1)F^{\omega, \omega'}_{h-1}(\lambda+1)$ has no pole in $-\alpha - \mathbb{N}$ for any $h \in \mathbb{Z}$ and any $\omega' \in \Omega^{n+1}_0$.$\hfill \blacksquare$\\
\end{enumerate}
\end{prop}

An easy consequence of the proposition above is the following:

\begin{cor}\label{pas de pole}
Under the hypothesis $H(\alpha,1)$ assume that the meromorphic extension of the holomorphic function  $F^{\omega, \omega'}_{h}(\lambda)$ has never a pole of order $\geq p$ at each point in $-\alpha - \mathbb{N}$  for some given $\omega' \in \Omega^{n+1}_0$  but for each $h \in \mathbb{Z}$. Then the same is true replacing $\omega$ by  any $w \in \Omega_0^{n+1}$ such that $[w]$ is in the fresco  $\mathcal{F}_\omega = B[a] [\omega] \subset H^{n+1}_0$.\\
\end{cor}

\parag{Proof of Corollary \ref{pas de pole}} Assume that the result is not true. So we have a $P \in B[a]$, an integer $m \in \mathbb{N}$ and some $h \in \mathbb{Z}$ such that $F^{P\omega, \omega'}_h(\lambda)$ has a pole of order at least equal to $p$ at the point $-\alpha -m$. First remark that if $ p+q \geq m+1$ the points $ii)$ and $iii)$ of Proposition \ref{tool 1} show that $F^{a^pb^q\omega, \omega'}_h(\lambda)$ has no pole at the point in $-\alpha - m$. So we may assume that the total degree of $P$ in $(a,b)$ is bounded by $m+1$ and the previous proposition gives a contradiction with our assumption.$\hfill \blacksquare$\\

 The  following important tool for the sequel is also a consequence of Proposition \ref{tool 1}, using  the Structure Theorem  for frescos of \cite{[B.09]}  extended  to the convergent case (see point 2 in the beginning of Section 3.1). \\
 
 The following terminology will be convenient:
 
 \parag{The property $P(\omega, \omega', p)$} In the standard  situation with the hypothesis $H(\alpha, 1)$ fix  two holomorphic germs $\omega, \omega'$ in $\Omega^{n+1}_0$. Let $p \geq 1$ be an integer and assume that there exists $h \in \mathbb{Z}$ such that $F^{\omega,\omega'}_h(\lambda)$ has a pole of order at least equal to  $p$ at a point in $-\alpha - \mathbb{N}$. Then we shall say that the integer  $m$ has the property $P(\omega, \omega', p)$ when $m$ is the smallest integer such that there exists an integer $h \in \mathbb{Z}$ with  a pole of order $\geq p$ at the point $\lambda = -\alpha-m$ for $F^{\omega, \omega'}_h(\lambda)$.

   \begin{prop}\label{tool 1bis}
   In the situation described above, assume that, for some $h \in \mathbb{Z}$, there exists a pole of order $\geq p$ at the point $-\alpha - m$ for $F^{\omega,\omega'}_h(\lambda)$. 
  Then the following properties hold true:
   \begin{enumerate}
   \item  Assume   that the integer  $m$  satisfies the property $P(\omega, \omega', p)$. For each $S \in \widehat{B}$ such that $S(0)\not= 0$ there exists a pole of order at least equal to  $p$ for $F^{S(b)\omega,\omega'}_{h+1}(\lambda)$ at the point $-\alpha - m$. Moreover, the integer $m$  satisfies  also the property $P(S\omega, \omega', p)$.
   \item  If $\mu \not= \alpha+m$   there exists a pole of order at least equal to  $p$ for $F^{(a -\mu b)\omega,\omega'}_{h+1}(\lambda)$ at the point $-\alpha - m-1$. Moreover, if  the integer  $m$  satisfies the property $P(\omega, \omega', p)$,  the integer $m+1$ satisfies the property $P((a - \mu b)\omega, \omega', p)$.
 \item For $\mu = \alpha+m$,  there exists a pole of order at least equal to  $p-1$ for $F^{(a -\mu b)\omega,\omega'}_{h+1}(\lambda)$ at the  point $-\alpha - m-1$. 
\end{enumerate}
     \end{prop}

     \parag{Proof} Assume that $X$ is a sufficiently small open neighborhood of \ $0$ \ in  \ $\mathbb{C}^{n+1}$ such that the germs $\omega$ and $\omega'$ are holomorphic on $X$ and that there exists  $ u \in \Omega^n(X)$ satisfying $du = \omega$ on $X$.\\
     Thanks to Stokes Formula and hypothesis $H(\alpha,1)$  (see Proposition 3.1.1 in  \cite{[B.22]} or Proposition    \ref{tool 1} above) the meromorphic function
     $$ F^{b\omega, \omega'}_h(\lambda) + (\lambda+1)F^{\omega, \omega'}_{h-1}(\lambda+1) = - \frac{1}{\Gamma(\lambda)}\int_X \vert f\vert^{2\lambda} \bar f^{-h} d\rho\wedge u\wedge \bar \omega' $$
     has no poles at points in $-\alpha -\mathbb{N}$ for any choice of $\omega', h$ and $\rho \in \mathscr{C}^\infty_c(X)$ which is identically $1$ near the origin (since $0$ is not in the support of $d\rho$).\\
     Since $m$ satisfies Property $P(\omega, \omega', p)$,  it is clear that for any positive integer $q$,  $F^{b^q\omega, \omega'}_{h'}(\lambda) $ has no pole of order $\geq p$ at $-\alpha -m'$ for each $m' \leq m-q$. Since we have never a  pole for $F^{\omega,\omega'}_h(\lambda)$  at points where 
     $\Re(\lambda) \ge 0$, we conclude that for any $S \in \widehat{B}$ with $S(0) \not= 0$ we have  a pole of order $p$ for $F^{S(b)\omega,\omega'}_{h+1}(\lambda)$ at the point $-\alpha -m$. Moreover    $m$ satisfies the property $P(S\omega, \omega', p)$. 
  
     With the same arguments (and the same Proposition 3.1.1 in  \cite{[B.22]} or Proposition \ref{tool 1}) the meromorphic function
     $$ F^{(a - \mu b)\omega, \omega'}_h(\lambda) - (\lambda + \mu + 1) F_{h-1}^{\omega, \omega'}(\lambda+1) $$
     has no pole at points in $-\alpha -\mathbb{N}$ for any choice of $\omega', h$ and $\rho \in \mathscr{C}^\infty_c(X)$ which is identically $0$ near the origin.\\
     Now the same line of proof gives the assertions 2 and 3 of the proposition using point $iv)$ in Proposition \ref{tool 1}.$\hfill \blacksquare$\\
     
     Here appears the main strategy of proof to locate the bigger order $p$ pole in $-\alpha-\mathbb{N}$ for a given pair $\omega, \omega'$.
          
     \begin{cor}\label{tool 2}
      Assume that there exists a pole of order at least equal to  $p$ at the point $-\alpha -m$ for $F^{\omega, \omega'}_h(\lambda)$  for some integer $h \in \mathbb{Z}$ and assume that the integer $m$ satisfies Property $P(\omega, \omega', p)$. \\
       Let $\Pi := (a -\mu_1b)S_1(a-\mu_2b)S_2 \dots (a-\mu_kb)S_k$  where $S_1, \dots, S_k$ are invertible elements in $B$ and $\mu_1, \dots, \mu_k$ are positive rational numbers.
      \begin{enumerate}
       \item Assume that  $\mu_j + j - k \not= \alpha +m$ for each $j \in [1,k]$.
       Then $F_{h+k}^{\Pi\omega, \omega'}(\lambda)$ has a pole of order at least equal to $p$ at the point $-\alpha - m - k$. Moreover  the integer 
       $m + k$ satisfies the property $P(\Pi\omega, \omega', p)$ .
       \item   If  $\mu_1$ is the only value of $j \in [1, k]$ such that $\mu + j - k = \alpha+m$ then $F_{h+k}^{\Pi\omega, \omega'}(\lambda)$ has a pole of order  $p-1$ at the point $-\alpha - m - k$.
      \end{enumerate}
      \end{cor}

    \parag{Proof} Using inductively the assertions $1, 2, 3$ of the previous proposition gives this corollary.$\hfill \blacksquare$\\
    
    \begin{cor}\label{tool 3}
    Assume $H(\alpha, 1)$ and that the nilpotent order of $(B[a] \omega)^{[\alpha]}$ (the $[\alpha]$-primitive quotient\footnote{The nilpotent order of $E^{[\alpha]}$ and of $E_{[\alpha]}$ are the same since the natural injection $E_{[\alpha]} \to E^{[\alpha]}$ has a finite co-dimensional co-kernel, since they have same rank; see Lemma 6.3.6 in \cite{[part I]}).} of the fresco $B[a] \omega$)  is at most $p-1$. Then for any choice of $\omega'$ and $ h$,  the meromorphic extension of $F^{\omega,\omega'}_h(\lambda)$ has no pole of order $\geq p$ at each  point  in  $-\alpha - \mathbb{N}$.
    \end{cor}

    \parag{Proof}  We shall prove the result by induction on $p \geq 1$. For $p = 1 $ our hypothesis means that    $(B[a]\omega)_{[\alpha]} = \{0\}$ so if $ \Pi := (a -\mu_1b)S_1(a-\mu_2b)S_2 \dots (a-\mu_kb)S_k$  where $S_1, \dots, S_k$ are invertible elements in $B$, generates the annihilator of $[\omega]$ in the geometric (a,b)-module $H^{n+1}_0$,  we may assume  that $\mu_1, \dots, \mu_k$ are not in $-\alpha -\mathbb{N}$. Then, since $F^{\Pi\omega, \omega'}_h(\lambda)$ has no poles  in $-\alpha - \mathbb{N}$ (see Proposition \ref{tool 1} (i)),  we obtain immediately a contradiction with  the assertion of Corollary \ref{tool 2} if we assume that for some choice of $\omega'$ and $h$ the meromorphic function has a pole at some point $-\alpha - m$.\\
    Thanks to the case proved above, we may replace $\omega$ by a generator of the fresco  $B[a] \omega^{[\alpha]}$, which means that we may assume now that $B[a]\omega$ is an $[\alpha]$-primitive fresco with nilpotent order at most $p-1$ with $p \geq 2$ (see again the previous footnote 11).\\
    Define $\mathcal{F} := S_{p-1}(B[a] \omega)\big/S_{p-2}(B[a] \omega)$. This fresco is $[\alpha]$-primitive,  semi-simple and generated by $[\omega]$. So the generator $\Pi := (a -\mu_1b)S_1(a-\mu_2b)S_2 \dots (a-\mu_kb)S_k$  where $S_1, \dots, S_k$ are invertible elements in $B$, of the  annihilator of the class $[\omega]$ in this semi-simple fresco\footnote{Recall that the saturation of this semi-simple fresco is a direct sum of rank $1$ geometric (a,b)-modules corresponding to pairwise distinct numbers since the Bernstein polynomial of  a fresco $\F$ is the characteristic polynomial of $-b^{-1}a$ acting of $\F^\sharp\big/b\F^\sharp$.} may be chosen  such that  we have any order for the sequence $\mu_j +j$. Since these numbers are pairwise distinct  there exists at most one $j \in [1, k]$ such that $\mu_j+j - k = \alpha + m $.  So we have have only  two possibilities: \\
     either there is no such $j \in [1, k]$ or there exists a unique $j \in [1, k]$ such that  $\mu_{j_0} + j_0 -k= \alpha + m $ and in this case we may choose $j_0 = 1$. \\
    So  using inductively Corollary \ref{tool 2} we see that if we assume that  $F^{\omega, \omega'}_h(\lambda)$ has a pole of order $\geq p$ at the point $-\alpha -m$, we shall find a pole of order $\geq p-1$ for $F^{\Pi\omega,\omega'}_{h+k}(\lambda)$ at the point $-\alpha - m - k$. Since the fresco $\mathcal{G}$ generated by the class $\Pi[\omega]$  satisfies $\mathcal{G} = S_{p-2}(B[a] [\omega])$, its nilpotent order is at most equal to $p-2$.  This contradicts our induction hypothesis.\\
    The case where there is no $j \in [1, k]$ such that $\mu_j + j - k = \alpha + m$ leads to a pole of order $\geq p$ at the point $-\alpha - m - k$, so gives also a contradiction.    $\hfill \blacksquare$\\
    
  \subsection{The final key}

Note that in Section 4.2 we always assume the existence of poles at some point in $-\alpha - \mathbb{N}$ for $F^{\omega, \omega'}_h(\lambda)$ (under our hypothesis $H(\alpha, 1)$)  and obtain consequences on the Bernstein polynomial of the fresco $\mathcal{F}_\omega$. These results go in the same direction than the results in \cite{[B.22]}. To go in the other direction, that is to say to prove the existence of such poles as consequence of informations on the Bernstein polynomial of $\mathcal{F}_\omega$, we shall use now the main idea of \cite{[B.84-a]} (and  also \cite{[B.84-b]}  in the case $\alpha = 1$). This is the point where the use of convergent (a,b)-modules is essential. It allows to show that the non vanishing of the class induced by $\omega$ in the $[\alpha]$-primitive quotient of $H^{n+1}_0$ implies that the cohomology class induces by $\omega/df$ in the spectral part for the eigenvalue $\exp(2i\pi\alpha)$ of the monodromy of $f$ acting on $H^n(F_0, \mathbb{C})$   does not vanish, where $F_0$ is the Milnor's fiber of $F$ at $0$.

\begin{thm}\label{effective 1}
Assume that $H(\alpha,1)$ is satisfied by $f : X \to D$, a Milnor representative  of a holomorphic germ near the origin in $\mathbb{C}^{n+1}$. Let $u \in \Omega^n(X)$ such that  there exists $m \in \mathbb{N}$ with $fdu = (\alpha+m)df\wedge u$ on $X$ and assume that the class induced by $u$ in $H^n(F_0, \mathbb{C})$ is not $0$. Then there exists a germ $\omega' \in \Omega^{n+1}_0$ and an  integer $h \in \mathbb{N}$ such that for any $\rho \in \mathscr{C}_c^\infty(X)$ which is identically $1$ near $0$ and with support small enough in order that $\rho\omega'$ is in $\mathscr{C}_c^\infty(X)^{n+1, 0}$, the meromorphic extension of
\begin{equation}
 \frac{1}{\Gamma(\lambda)} \int_X \vert f\vert^{2\lambda}\bar f^{-h}\rho \frac{df}{f}\wedge u \wedge \bar{\omega'} 
 \end{equation}
has a pole at $-\alpha -m$.
\end{thm}

\parag{Proof}  Define, for $j \in \mathbb{N}$, the $(n, 0)-$current on $X$ by the formula\footnote{Here $Pf(\lambda = \lambda_0, F(\lambda))$ denotes the constant term in the Laurent expansion at $\lambda = \lambda_0$ of the meromorphic function $F(\lambda)$.}
$$ \langle T_j, \psi \rangle := Pf\big(\lambda = -\alpha -m, \  \frac{1}{\Gamma(\lambda)} \int_X \vert f\vert^{2\lambda}\bar f^{-j}u\wedge \psi \big)$$
where $\psi$ is a test form of type $(1,n+1)$ which is $\mathscr{C}_c^\infty$ in $X$. 
\parag{Claim} Then we have the following properties for each $j \in \mathbb{N}$
\begin{enumerate}
\item  $\bar{f}T_{j+1} = T_j$ on $X$
\item The support of the current $d'T_j $ is contained in $\{0\}$.
\item The support of the current $d''T_j + (\alpha+m+j)\bar{df}\wedge T_{j+1}$ is contained in $\{0\}$.
\end{enumerate}
\parag{proof of the claim} The first assertion is clear.\\ 
Let us compute $d'T_j$. Let $\varphi$ be a $\mathscr{C}^\infty_c(X)$ test form of type $(0,n+1)$. We have
\begin{equation*}
 \langle d'T_j, \varphi \rangle := (-1)^n\langle T_j, d'\varphi \rangle = (-1)^n\langle T_j, d\varphi \rangle
 \end{equation*}
 But for $\Re(\lambda) \gg 1$ the form $\vert f\vert^{2\lambda}\bar f^{-j.}u\wedge \varphi$ is in $\mathscr{C}^1_c(X)$  and Stokes Formula and the meromorphic continuation give
 $$ 0 = \frac{1}{\Gamma(\lambda)} \int_X d\big(\vert f\vert^{2\lambda}\bar f^{-j}u\wedge \varphi\big) = \frac{(\lambda+\alpha +m)}{\Gamma(\lambda)} \int_X \vert f\vert^{2\lambda}\bar f^{-j} \frac{df}{f}\wedge u\wedge \varphi + (-1)^n\langle T_j, d'\varphi \rangle $$
 because $du = (\alpha+m)\frac{df}{f}\wedge u$ and $d\bar{f}\wedge \varphi \equiv 0$. Then we obtain\footnote{Here $Res(\lambda = \lambda_0, F(\lambda))$ denotes the residue at $\lambda = \lambda_0$ of the meromorphic function $F(\lambda)$.}
 \begin{equation}
\langle  d'T_j, \varphi \rangle =  Res\Big(\lambda = -\alpha-m,\  \frac{1}{\Gamma(\lambda)} \int_X \vert f\vert^{2\lambda}\bar f^{-j}\frac{df}{f}\wedge u\wedge \varphi \Big)
\end{equation}
This gives our assertion $2$ because we know that the poles of the meromorphic extension of $\frac{1}{\Gamma(\lambda)} \int_X \vert f\vert^{2\lambda}\bar f^{-j}\square$ at points in $-\alpha + \mathbb{Z}$ are supported by the origin, thanks to our hypothesis $H(\alpha,1)$.\\
In an analogous way let us compute $d''T_j$; let $\psi$ be a $\mathscr{C}^\infty_c(X)$ test form of type $(1, n)$. We have:
\begin{equation*}
 \langle d''T_j, \psi \rangle := (-1)^n\langle T_j, d''\psi \rangle = (-1)^n\langle T_j, d\psi\rangle
 \end{equation*}
But for $\Re(\lambda) \gg 1$ the form $\vert f\vert^{2\lambda}\bar f^{-j}u\wedge \psi$ is in $\mathscr{C}^1_c(X)$ and so:
 $$d\big(\vert f\vert^{2\lambda}\bar f^{-j}u\wedge \psi\big) = (\lambda - j)\vert f\vert^{2\lambda}\bar f^{-j-1}\bar{df}\wedge u\wedge \psi + (-1)^n\vert f\vert^{2\lambda}\bar f^{-j}u\wedge d\psi $$
 because the type of $du$ as well as the type of $df\wedge u$ is $(0, n+1)$. Then Stokes Formula and the meromorphic continuation give
 $$ \langle d''T_j + (\alpha+m+j)\bar{df}\wedge T_{j+1}, \psi \rangle = Res\Big(\lambda = -\alpha-m, \ \frac{1}{\Gamma(\lambda)} \int_X \vert f\vert^{2\lambda}\bar f^{-j-1}d\bar{f}\wedge u\wedge \psi \Big) .$$
 This proves the assertion $3$, again thanks to our hypothesis $H(\alpha,1)$.\\
 Now we shall argue by contradiction and we shall assume that for each  $j_0 \in \mathbb{N}$,  the current $d'T_{j_0}$ induces the class zero in the conjugate of the space $H^{n+1}_{[0]}(X, \mathcal{O}_X)$ which means that there exists  a $(n, 0)-$current $\Theta_{j_0}$ with support $\{0\}$ satisfying $d'\Theta_{j_0} = d'T_{j_0}$ on $X$. Then, as we have $\bar f^kT_{j_0} = T_{j_0-k}$ for any $k \in \mathbb{N}$ thanks to $1$, we obtain that $d'\bar f^k\Theta_{j_0} = \bar f^kd'\Theta_{j_0} = \bar f^k T_{j_0} = T_{j_0-k}$. Now we fix $j_0 \gg 1$ and   define, for each $j  \leq j_0$, $\Theta_j :=\bar f^{j_0 -j} \Theta_{j_0}$. 
So for any such $j \leq j_0$ this gives $d'\Theta_j = d'T_j$.\\
 Now we shall use  Lemma $C_1, C_2$  and Lemma $D$ in \cite{[B.84-a]}  and Lemma $C'_1, C'_2$  in  \cite{[B.84-b]} in the case $\alpha = 1$, for the family of currents $\tilde{T}_{j} := T_j - \Theta_j$ for $j \leq j_0$ to obtain a contradiction.\\
 These currents  satisfy
 \begin{enumerate}
 \item $d'\tilde{T}_j = 0$ on $X$.
 \item  $d''\tilde{T}_j + (\alpha + m + j) \bar{df}\wedge \tilde{T}_{j+1}$ has its support in $\{0\}$.
 \item  The current $\tilde{T}_j$ coincides with $\vert f\vert^{-2(\alpha +m)}\bar{f}^{-j} u$ on the Milnor fiber $F_0 = f^{-1}(s_0)$ (these currents are smooth outside $Y$).
 \end{enumerate}
 Note that we have $H^p(X \setminus \{0\}, \mathcal{O}_X) = 0$ for $1 \leq p \leq n-1$ which is used for  checking the hypothesis of Lemmas $C'_1, C'_2$ in the case $\alpha = 1$ .\\
 Then the cited Lemma  contradict our assumption  that the class induced by $u$ in $H^n(F_0, \mathbb{C})$  does not vanish.\\
 So we obtain that there exists $j_0 \in \mathbb{N}$ such that  the class induced by $d'T_{j_0}$  does not vanish in the dual of the space $\overline{\Omega^{n+1}_0}$ of the germs at the origin of anti-holomorphic volume forms on $\mathbb{C}^{n+1}$ (and then  for any $j_0 +k$ for $k \in \mathbb{N}$ also). So there exists $\omega' \in \Omega^{n+1}_0$ and $\rho \in  \mathscr{C}_c^\infty(X)$ which is identically $1$ near $0$ and with support small enough in order that $\rho\omega'$ is in $\mathscr{C}_c^\infty(X)^{n+1,0}$ such that
 $$ \langle d'T_{j_0}, \rho\overline{\omega}' \rangle = Res\Big(\lambda = -\alpha-m, \ \frac{1}{\Gamma(\lambda)} \int_X \vert f\vert^{2\lambda}\bar f^{-j_0}\frac{df}{f}\wedge u\wedge \rho \overline{\omega}' \Big) \not= 0 $$
 concluding the proof of the theorem.$\hfill \blacksquare$
 
 \parag{Remark} To apply the results in degree $n$  of \cite{[B.84-a]} and \cite{[B.84-b]} used in the previous proof, it is enough to assume that $j_0 = n+1$ to conclude. That is to say that if there is no pole in the range $ [-\alpha -n-1, -\alpha] \cap \{-\alpha - \mathbb{N}\} $ for $ F^{\omega, \omega'}_h(\lambda)$, there is no pole in $-\alpha - \mathbb{N}$. 
 
\begin{cor}\label{effective 2}
Assume that we have holomorphic forms $u_j \in  \Omega^n(X)$ for each integer  $j$   in $[-N,p]$ such that
$$ fdu_j = (\alpha+m)df\wedge u_j +  df\wedge u_{j-1}  $$ 
with the hypothesis that $[du_j] = 0$ in $H^{n+1}_0$ for each $j \in [-N, 0]$ and  such that $u_1 $ induces   in $H^n(F_0, \mathbb{C})$ a class which  is not $0$.\\
 Then there exists $h \in \mathbb{N}$ and $\omega' \in \Omega^{n+1}_0$ such that the meromorphic extension of
\begin{equation}
 \frac{1}{\Gamma(\lambda)} \int_X \vert f\vert^{2\lambda}\bar f^{-h}\rho \frac{df}{f}\wedge u_p \wedge \bar{\omega'} 
 \end{equation}
has a pole of order at least equal to $p$  at  the point $\lambda = -\alpha -m$.
\end{cor}

\parag{Proof} For $\Re(\lambda) \gg 1$ the differential form $\vert f\vert^{2\lambda}\bar f^{-h}\rho u_j\wedge\bar \omega'$ is of class $\mathscr{C}^1$ and satisfies
\begin{align*}
&  d(\vert f\vert^{2\lambda}\bar f^{-h}\rho u_j\wedge\bar \omega') = (\lambda+\alpha+m)\vert f\vert^{2\lambda}\bar f^{-h}\rho (df/f)\wedge u_j\wedge \bar \omega' + \\
&  + \vert f\vert^{2\lambda}\bar f^{-h}\rho (df/f)\wedge u_{j-1}\wedge \bar \omega'   +   \vert f\vert^{2\lambda}\bar f^{-h}d\rho \wedge u_j\wedge \bar \omega' .
\end{align*}
Then Stokes Formula and the meromorphic extension gives, where $P_q(\lambda = \lambda_0, F(\lambda))$ means the coefficient of $(\lambda-\lambda_0)^{-q}$ in the Laurent expansion\ of the meromorphic function $F$ at the point $\lambda = \lambda_0$, that for each $q \geq 1$ and each $j$  we have:
\begin{align*}
& P_{q+1}\Big((\lambda = -\alpha-m, \frac{1}{\Gamma(\lambda)} \int_X \vert f\vert^{2\lambda}\bar f^{-h}\rho (df/f)\wedge u_j\wedge \bar \omega')\Big) =  \\
& - P_q\Big(\lambda = -\alpha-m, \frac{1}{\Gamma(\lambda)} \int_X \vert f\vert^{2\lambda}\bar f^{-h}\rho (df/f)\wedge u_{j-1}\wedge \bar \omega'\Big).
\end{align*}
Then the fact that there exists $h \in \mathbb{N}$ and $\omega' \in \Omega^{n+1}_0$ with (here $P_1 = Res$ !)
$$P_1\Big(\lambda = -\alpha-m, \frac{1}{\Gamma(\lambda)} \int_X \vert f\vert^{2\lambda}\bar f^{-h}\rho (df/f)\wedge u_{1}\wedge \bar \omega'\Big) \not= 0 $$
completes the proof of the theorem using the formulas above with $q \in [1, p-1]$ with $j = q+1$.$\hfill \blacksquare$\\

To be able to use the previous corollary, the following lemma, combined with Corollary \ref{Bloc de Jordan dans fresque} will be useful.

\begin{lemma}\label{Bloc Jordan formes} Let $w_1, \dots, w_p$ be in $\Omega^{n+1}_0$ such that the induced class in $H^{n+1}_0$   satisfy the relations:
\begin{equation*}
	a[w_j] = (\alpha +m)b[w_j] + b[w_{j-1}] \quad \forall j \in [1,p] \quad {\rm with \ the \ convention} \quad [w_{0}] = 0 \tag{$\star$}
\end{equation*}
 Then there exists an integer $N$ and  $u_1 \dots, u_p$ in $\Omega^n_0$ such that
\begin{equation*}
fdu_j = (\alpha+m +N)df\wedge u_j + df\wedge u_{j-1}   \quad {\rm with \ the \ convention} \quad u_0 = 0 \tag{$\star\star$} 
\end{equation*}
and  such that we have \ $[du_j] = (a+b)^{N}[w_j]$ \  in $H^{n+1}_0$.
\end{lemma}

\parag{Proof} Choose for each $j \in [1,p]$ a  $v_j \in \Omega^n_0$ such that $dv_j = w_j$. Then for each $j \in [1, p]$ the class induced in $H^{n+1}_0$ by the form
$$fdv_j - (\alpha + m)df \wedge v_j - df \wedge v_{j-1} \quad {\rm with \ the \ convention}  \ v_0 = 0 $$
is of $a$-torsion in  $\mathcal{H}^{n+1}_0$. So there exists an integer $N$  and $t_j \in (Ker\, df)^n_0$ such that, for $j \in [1,p]$, we have
$$ f^{N+1}dv_j - (\alpha + m)df \wedge f^Nv_j - df \wedge f^N v_{j-1} = fdt_j .$$
This equality may be written
$$ fd(f^Nv_j + t_j) -(\alpha + m + N) df\wedge(f^Nv_j + t_j) - df\wedge (f^Nv_{j-1} + t_{j-1}) = 0 $$
with the convention $t_0 = 0$, using the fact that $df\wedge t_j = 0$ for each $j$. Then defining $u_j := f^Nv_j + t_j$ for $j \in [1, p]$ concludes the proof since the class induced in $H^{n+1}_0$ by
$du_j$ is equal to $a^N[w_j] + N a^{N-1}b[w_j]$, thanks to the equality $a^N + Na^{N-1}b = (a+b)^N$ (see the exercise below).$\hfill \blacksquare$
 
\parag{Exercise} Show that the commutation relation $ab - ba = b^2$  implies the relation
$$ (a+b)^q = a^{q-1}(a + qb) \quad \forall q \in \mathbb{N}^*.$$
\smallskip

Note that $b(a+b)^N = a^Nb$ implies that in a geometric (a,b)-module an element  $x$ which satisfies $(a+b)^Nx= 0$ is nul, since $a$ and $b$ are injective.

\subsection{Statements and proofs}

Our first result gives an improvement of the result in \cite{[B.22]} but is also a precise converse of this statement. It shows the interest in considering the higher  Bernstein polynomials introduced in section 2.

The content of the following result is the direct part of Theorem \ref{1} in the introduction.

\begin{thm}\label{Higher Bernstein 1} In the standard situation described above, assume that the hypothesis  $H(\alpha,1)$ is satisfied. Consider  a germ $\omega \in \Omega^{n+1}_0$ such that  the $p$-th Bernstein polynomial of the fresco  $\mathcal{F}_\omega := B[a] \omega$ in $ H^{n+1}_0$ has a root in $-\alpha -\mathbb{N}$. Then there exists $\omega' \in \Omega^{n+1}_0$ and an integer $h \in \mathbb{Z}$ such that the meromorphic extension of the integral
\begin{equation*}
F^{\omega'}_{\omega,h}(\lambda) :=  \frac{1}{\Gamma(\lambda)}\int_X \vert f\vert^{2\lambda} \bar f^{-h} \rho \omega\wedge\bar \omega' \tag{A}
\end{equation*}
has a pole of order at least equal to $p$ at $\lambda = -\alpha - m$ for $m$ a large enough integer, where $\rho \in \mathscr{C}^\infty_c(X)$ is identically $1$ near zero.\\
\end{thm}

\parag{Remark}  The converse of this result, so the second part of Theorem \ref{1} in the introduction, so the fact  that,  for a germ $\omega \in \Omega^{n+1}_0$, the existence of such \  $\omega', h, m$ \  giving a pole of order $p$ at a point in $-\alpha - \mathbb{N}$ for $(A)$  implies that the $p$-th Bernstein polynomial of the fresco $\mathcal{F}_\omega = B[a] \omega$ has a root in $-\alpha -\mathbb{N}$,  will be a consequence of the Theorem \ref{sauts} which is more precise, using the following result proved in Corollary \ref{23/12}    
\begin{itemize}
\item If the $q$-th Bernstein polynomial of the geometric (a,b)-module  $\E$ has a root in $-\alpha - \mathbb{N}$ then for each $p \in [1,q]$ the $p$-th Bernstein polynomial of $\E$ has also a root in  $-\alpha - \mathbb{N}$.\\
\end{itemize}

For the proof of Theorem  \ref{Higher Bernstein 1} we  shall need the following result.

\begin{prop}\label{crucial}
Assume that the hypothesis $H(\alpha, 1)$. Suppose that $u_1 \in \Omega^n_0$ satisfies the relation $fdu_1 = (\alpha + m) df \wedge u_1$ for some integer $m$. If $[du_1]$ is not zero in $H^{n+1}_0$ then the cohomology class induced by $u_1$ in $H^n(F_0, \mathbb{C})$ is not zero. So $u_{\vert F_0}$  induces a class which  is an eigenvector of the monodromy for the eigenvalue $\exp(-2i\pi\alpha)$.
\end{prop}

\parag{Proof} Thanks to Grothendieck (see \cite{[G.65]}), the meromorphic relative  de Rham  complex of $f$ computes the cohomology of $X \setminus f^{-1}(0)$ and under the hypothesis $H(\alpha, 1)$ the spectral sub-space $H^n(F_0,\mathbb{C})_{\exp(-2i\pi\alpha)}$ of the monodromy is isomorphic to the $n$-th cohomology group of the complex 
$$\left( \Omega_0^{\bullet}[ f^{-1}]\big/df\wedge \Omega_0^{\bullet -1}[ f^{-1}], (d - \alpha \frac{df}{f} \wedge)^{\bullet}\right).$$
If we assume that $u_1$ induces $0$ in $H^n(F_0, \mathbb{C})$, since we have
 $$d(f^{-m}u_1) - \alpha \frac{df}{f} \wedge f^{-m}u_1 = 0,$$
 there exists $v, w \in \Omega^{n-1}_0[ f^{-1}]$ such that
$$ dv - \alpha \frac{df}{f} \wedge v = f^{-m}u_1 + df \wedge w .$$
This gives 
\begin{align*}
& d(f^{-m}u_1) = -df \wedge dw + \alpha df\wedge v/f \quad {\rm and \ then} \\
& f^{-m}du_1 - m \frac{df}{f} \wedge f^{-m}u_1 = (1 - m/(m +\alpha)) f^{-m}du_1 = df \wedge d(-w + \alpha v/f) 
\end{align*}
and this implies, since $\alpha$ is in $]0, 1]$,  that $[du_1]$ is of $a$-torsion in $\mathcal{H}^{n+1}_0$ and then  $0$ in $H^{n+1}_0$. Contradiction.$\hfill \blacksquare$\\

\parag{Proof of Theorem \ref{Higher Bernstein 1}} Using Corollary \ref{Bloc de Jordan dans fresque} there exist $[w_1], \dots, [w_p]$ in $\mathcal{F}_\omega$  and an integer $m \in \mathbb{N}$ satisfying the relations:
\begin{equation*}
	a[w_j ]= (\alpha +m)b[w_j ]+ b[w_{j-1}] \quad \forall j \in [1,p] \quad {\rm with \ the \ convention} \quad [w_0] = 0 \tag{*}
\end{equation*}
and which are $B$-linearly independent in $\mathcal{F}_\omega$. Assuming that the Theorem does not hold would imply, thanks to Corollary \ref{pas de pole} and to Corollary \ref{effective 2}, that writing $w_1 = du_1$ with $u_1 \in \Omega^n_0$,  the class induced by $u_1$ in $H^n(F_0, \mathbb{C})$ vanishes. \\
But this contradicts the hypothesis that $[w_1] $ is not  zero  in
 $\mathcal{F}_\omega \subset H^{n+1}_0$, thanks to Proposition \ref{crucial}.$\hfill \blacksquare$\\

The following corollary of Theorem  \ref{Higher Bernstein 1} is clear since we may use a  Bernstein identity at the origin to describe the poles of the meromorphic extension of the distribution $\frac{1}{\Gamma(\lambda)}\vert f \vert^{2\lambda} \bar f^{-h}$ for any $h \in \mathbb{Z}$ (see \cite{[B.81]} or \cite{[Bj.93]}).

\begin{cor}\label{Higher Bernstein 2} In the situation of the previous theorem, the existence of a germ $\omega \in \Omega^{n+1}_0$ such that  the $p$-th Bernstein polynomial of the fresco  $\Abc \omega \subset H^{n+1}_0$ has a root in $-\alpha -\mathbb{N}$ implies the existence of at least $p$ roots of the reduced $b$-function $b_{f,0}$ of $f$ at the origin in $-\alpha -\mathbb{N}$ counting multiplicities.$\hfill \blacksquare$\\
\end{cor}

\parag{Remark} The interest of this corollary lies in the fact that the existence of $p$ roots in $-\alpha - \mathbb{N}$ \  for the reduced Bernstein polynomial $b_{f,0}$ does not implies, in general under our hypothesis, the existence of a pole of order $p$ at some point  $\lambda = -\alpha - m$ with  $m \in \mathbb{N}$ \ large, for the meromorphic extension of $ \frac{1}{\Gamma(\lambda)}\int_X \vert f\vert^{2\lambda} \bar f^{-h}\varphi$ for some test $(n+1,n+1)$-form $\varphi$. \\
The consideration of higher  Bernstein polynomials of frescos associated to germs $\omega \in \Omega^{n+1}_0$ is then  a tool which may help to determine the nilpotent order of the monodromy of $f$ at the origin in the case of an isolated singularity for the eigenvalue $\exp(2i\pi\alpha)$. \\

Our next result is an improvement of Theorem \ref{Higher Bernstein 1}.

\begin{thm}\label{Higher Bernstein 2}
In the standard situation, we assume that the hypothesis $H(\alpha, 1)$ is satisfied.  Assume that there exists $\omega \in \Omega^{n+1}_0$ such that $B^p_{\mathcal{F}^{[\alpha]}}$ has a root  in $-\alpha - \mathbb{N}$ where $\F := B[a]\omega \subset H^{n+1}_0$ and where $p = d(\mathcal{F}^{[\alpha]})$ is the nilpotent order of the fresco  $\mathcal{F}^{[\alpha]}$. Let $-\alpha -m$ the biggest root of $B^p_{\mathcal{\mathcal{F}^{[\alpha]}}}$ in $-\alpha - \mathbb{N}$. Then there exists $ \omega' \in \Omega^{n+1}_0$ and  $h \in \mathbb{Z}$ such that $F^{\omega, \omega'}_h(\lambda)$ has a pole of order $p$ at the point  $-\alpha -m$.
\end{thm}

Recall that, of course, in the previous statement  $B^p_{\mathcal{F}}$ denotes the $p$-th Bernstein polynomial of the fresco $\mathcal{F}$.

\parag{Proof} First recall that, thanks to Lemma  6.3.6 in \cite{[part I]} ,  we have, for any geometric (a,b)-module $\mathcal{E}$, the equality
$d(\mathcal{E}_{[\alpha]}) = d(\mathcal{E}\big/\mathcal{E}_{[\not= \alpha]}) = d(\mathcal{E}^{[\alpha]})$.  \\
Let $-\alpha -m$ be the biggest root of $B^p_{\mathcal{F}^{[\alpha]}}$. Then we may choose a J-H. sequence of $\mathcal{F}^{[\alpha]} \big/S_{p-1}(\mathcal{F}^{[\alpha]} )$ such that its last quotient is isomorphic to $\mathcal{E}_{\alpha+m}$ . This possible because the fresco $\mathcal{F}^{[\alpha]} \big/S_{p-1}(\mathcal{F}^{[\alpha]} )$ is semi-simple  and has $-\alpha-m$ as a root of its Bernstein polynomial (remind  that the Bernstein polynomial fo $\mathcal{F}^{[\alpha]} \big/S_{p-1}(\mathcal{F}^{[\alpha]} )$ divides $B^p_{\mathcal{F}^{[\alpha]}}$ which also divides $B^p_{\mathcal{F}}$). Then if $\Pi_0$ is the generator of the annihilator of $[\omega]$ in $\mathcal{F}^{[\alpha]} \big/S_{p-1}(\mathcal{F}^{[\alpha]} )$, it may be written $\Pi_0 = (a - (\alpha+m+1-k)b) \Pi'_0$ where $k$ is the rank of $\mathcal{F}^{[\alpha]} \big/S_{p-1}(\mathcal{F}^{[\alpha]} )$. Then, choosing a J-H. sequence of $\mathcal{F}$ which begins by a J-H. sequence of $\mathcal{F}_{[\not=\alpha]}$ and  ending by the J-H. sequence of $\mathcal{F}^{[\alpha]} $ chosen above, we see that the annihilator of $\omega$ in $\mathcal{F}$ may be written as $\Pi = \Pi_2\Pi_1(a - (\alpha+m+1-k)b)\Pi'_0 $ with $B[a]\big/B[a] \Pi'_0 $  semi-simple $[\alpha]$-primitive with a Bernstein polynomial having roots strictly less than $-\alpha-m$, with $d(B[a]\big/B[a]\Pi_1) \leq p-1$, since this fresco is isomorphic to $S_{p-1}(\mathcal{F}^{[\alpha]})$ and with $(B[a]\big/B[a]\Pi_2)_{[\alpha]} = \{0\}$ since $B[a]\big/B[a]\Pi_2$ is isomorphic to $\mathcal{F}_{[\not=\alpha]}$.\\
Now, applying Theorem \ref{Higher Bernstein 1} we find $\omega' \in \Omega^{n+1}_0$, $h \in \mathbb{Z}$ and $m_1 \in \mathbb{N}$ such that $F^{\omega, \omega'}_h(\lambda)$ has an order $p$ pole at the point $-\alpha - m_1$ and such that  the integer $m_1$ satisfies the property $P(\omega, \omega', p)$. \\
Using then Proposition \ref{tool 1} we see that  if $m \not= m_1$ we obtain a  contradiction with  Corollary \ref{pas de pole} because we find a pole of order $p$ at a point $-\alpha - m_1 -k$ for the meromorphic extension of $F^{\Pi_0\omega, \omega'}_{h-q-1}(\lambda)$ where $k$ is the degree in $a$ of $\Pi_0$.\\
 So we obtain that  $m = m_1$, concluding the proof.$\hfill \blacksquare$ \\

The following corollaries are  obvious consequences of the previous result.

\begin{cor}\label{first pole 1}  In the standard situation described above, under the assumption $H(\alpha,1)$, consider a germ $\omega \in \Omega_0$ and assume that $-\alpha -m$ is the biggest possible pole in $-\alpha - \mathbb{N}$ for any choices of $\omega' \in \Omega^{n+1}_0$ and any $h \in \mathbb{Z}$ for the meromorphic functions $F^{\omega, \omega'}_h(\lambda)$. Then $-\alpha -m$ is the biggest root in $-\alpha - \mathbb{N}$ of the Bernstein polynomial of the fresco $\mathcal{F}_\omega := (B[a] \omega) \subset H^{n+1}_0$. $\hfill \blacksquare$ \\
\end{cor}

\begin{cor}\label{first pole 2}
In the standard situation described above, under the assumption $H(\alpha,1)$, assume that $-\alpha -m$ is the biggest root of the Bernstein polynomial in $-\alpha - \mathbb{N}$ of the geometric (a,b)-module $H^{n+1}_0$. Then there exists $h \in \mathbb{Z}$ such the meromorphic extension of the distribution $\frac{1}{\Gamma(\lambda)}\vert f\vert^{2\lambda} \bar f^{-h}$ has a pole at $-\alpha - m$.$\hfill \blacksquare$\\
\end{cor}

Note that it is enough to consider the integers $h \in [m+1, n-m]$ in the previous statement since the exponent of $\bar f$ has to be negative and thanks to the remark following Theorem \ref{effective 1} which  implies that $\alpha + m \leq n+1$.\\

The following consequence of the previous corollary is obvious, since in the case of an isolated singularity at $0$ for $f$ it is known (see \cite{[M.75]}) that the  Brieskorn module coincides with $H^{n+1}_0$ and that its Bernstein polynomial coincides with the reduced Bernstein polynomial  $\tilde{b}_f$ of $f$.

\begin{cor}\label{first pole 3}
Assume that  the germ $f : (\mathbb{C}^{n+1}, 0) \to (\mathbb{C}, 0)$ of holomorphic function  has an isolated singularity at the origin. For $\alpha \in ]0, 1] \cap \mathbb{Q}$ let $-\alpha -m$ be the biggest root of the reduced Bernstein polynomial of $f$ in $-\alpha -\mathbb{N}$. Then there exists $h \in \mathbb{Z}$ such the meromorphic extension of the distribution $\vert f\vert^{2\lambda} \bar f^{-h}\big/ \Gamma(\lambda) $ has a pole at $-\alpha - m$.$\hfill \blacksquare$
\end{cor}

\parag{Question} In the case of an isolated singularity for the eigenvalue $\exp(2i\pi\alpha)$ of the monodromy (so with  our hypothesis $H(\alpha, 1)$), is the Bernstein polynomial of the geometric (a,b)-module $ (H^{n+1}_0)^{[\alpha]} := H^{n+1}_0\big/(H^{n+1}_0)_{\not=[\alpha]}$ (which is the biggest polynomial having its root in $-\alpha - \mathbb{N}$ and  dividing the Bernstein polynomial of the (a,b)-module  $H^{n+1}_0$)  coincides with the biggest polynomial having its root in $-\alpha - \mathbb{N}$ and  dividing the reduced Bernstein polynomial of $f$ at the origin ?\\

  \subsection{Some improvements of Theorem 3.1.2 in \cite{[B.22]}}
   
   The goal of this paragraph is to show that, using the higher Bernstein polynomials of the fresco $\mathcal{F}_{f,\omega}$ generated by the class of $\omega$ in $H^{n+1}_0$, and the tools introduced above, we can improve the main result in \cite{[B.22]} Theorem 3.1.2. 
      
   We begin by some remarks to make clear the correspondence between our present notations with those used in \cite{[B.22]}.
   
   \parag{Remarks} \begin{enumerate}
   \item We use here the notation $H(\alpha, 1)$ with $\alpha \in ]0, 1] \cap \mathbb{Q}$ instead of the notation $H(\xi, 1)$ with $\xi \in \mathbb{Q}$.
   \item  To consider a form  $\psi \in \mathscr{C}^\infty_c(\mathbb{C}^{n+1})^{0, n+1}$  with small enough support and such that $d\psi = 0$ in a neighborhood of $0$ is equivalent to consider $\rho.\bar \omega'$ where $\omega' $ is in $\Omega^{n+1}_0$ and $\rho$ is a function in $\mathscr{C}^\infty_c(\mathbb{C}^{n+1})$ with small enough support which is identically $1$ near the origin. \\
Indeed any such $\psi$ may be written as $\psi = \bar \omega'$  for some $\omega' \in \Omega^{n+1}$ near the origin  thanks to Dolbeault' Lemma, and then $\psi -\rho \bar \omega'$ is identically $0$ near the origin, so replacing $\psi$ by $\rho \bar \omega' $ do not change the poles which may appear in $-\alpha - \mathbb{N}$  for the functions we are looking at (what ever is the choice of $h \in \mathbb{Z}$ thanks to our hypothesis $H(\alpha, 1)$). \\
Then, for $\omega, \omega'$ in $\Omega^{n+1}_0$  we use the notation $F^{\omega,\omega'}_h(\lambda)$ where the function $\rho \in {\mathscr{C}^\infty}_c(\mathbb{C}^{n+1})$ which is identically $1$ near the origin and has a sufficiently small support in order that $\rho \omega \wedge \bar \omega'$ is smooth, does not appear in this notation  because the poles at points in $-\alpha - \mathbb{N}$ do not depend on the choice of this $\rho$. \\
This corresponds to  the notation $F^\psi_h(\lambda)$ where $\psi$ is in $\mathscr{C}^\infty_c(X)^{0,n+1}$ is $d$-closed near the origin (where $\omega $ is given in $\Omega^{n+1}(X)$) and  with $\psi = \rho\bar \omega'$.
\item Note also that we change the sign of the integer $h \in \mathbb{Z}$ between these two articles.
\end{enumerate}

\begin{thm}\label{improve 1}
Let $\alpha \in ]0,1]$ and assume the hypothesis $H(\alpha, 1)$ for  the germ at the origin in $\mathbb{C}^{n+1}$ of holomorphic function $\tilde{f} : (\mathbb{C}^{n+1},0) \to (\mathbb{C}, 0)$. Assume that
 $\omega$ in $\Omega^{n+1}_0$ is such that there exists an integer $h \in \mathbb{Z}$ and a form $\omega' \in \Omega^{n+1}_0$ for which the function $F^{\omega, \omega'}_h(\lambda) $ has  a pole of order
  $p \geq 1 $ at some  point $\xi$  in $\{-\alpha -\mathbb{N}\}$. Note $\xi_p = -\alpha - m$ be the biggest such  number $\xi$ in  $\{-\alpha - \mathbb{N}\}$  for any choice of $\omega'$ and $h \in \mathbb{Z}$. Then the $p$-th Bernstein polynomial of the fresco $\mathcal{F}_{f, \omega} := B[a]\omega \subset H^{n+1}_0$ has a root in $[-\alpha -m, -\alpha] \cap \mathbb{Z}$.
\end{thm}

\parag{Proof} Note $P := P_1P_2$ the annihilator of the class of $[\omega]$ in the  $[\alpha]$-primitive quotient 
$$\F^{[\alpha]} := \F_{f,\omega}\big/\big(\F_{f,\omega}\big)_{\not= \alpha} $$
 of the fresco   $\F_{f, \omega}:= B[a][\omega]$ inside the (a,b)-module $H^{n+1}_0$ associated to $f$, where $P_2$ is the annihilator of $[\omega]$ in $\F^{[\alpha]}\big/S_{p-1}\big(\mathcal{F}^\alpha\big)$. If $F^{\omega,\omega'}_h(\lambda)$ has a pole of order at least equal to $p$ at the point $-\alpha -m$ and if $-\alpha-m$ is not a root of the $p$-th Bernstein polynomial of $\F^{[\alpha]}$, then $-\alpha -m$ is not a root of the (usual) Bernstein polynomial of the fresco $B[a]\big/B[a]P_2$ which is isomorphic to $S_{p-1}(\mathcal{F}^{[\alpha]})$. In this situation, using Corollary \ref{tool 2} we see  $F^{P_2\omega, \omega'}_{h+p_2}(\lambda)$ has a pole of order at least equal to $p$ at $-\alpha -m-k$, where $k$ is the rank of the fresco $\mathcal{F}^{[\alpha]}\big/S_{p-1}\big(\mathcal{F}^{[\alpha]}\big)$. But this is impossible, according to Corollary \ref{tool 3}  since the nilpotent order of  $S_{p-1}\big(\mathcal{F}^{[\alpha]}\big)$ is $p-1$. So $-\alpha-m$ is a root of some $(p+j)$-th Bernstein polynomial of $\mathcal{F}^{[\alpha]}$ for some integer j $\geq 0$.$\hfill \blacksquare$\\

The end of Theorem 3.1.2 in \cite{[B.22]} is also improved as follows:

\begin{cor}\label{improve 2} In the situation of the previous theorem, let, for  each integer $s$  in $[1,p]$,  $\xi_s$ be the biggest element in $-\alpha -\mathbb{N}$ for which there exists $h \in \mathbb{Z}$ and $\omega' \in \Omega^{n+1}_0$ such that $F^{\omega, \omega'}_h(\lambda)$ has a pole of order  at least equal to $s$ at $\xi_s$. Then $\xi_s$ is a root of  some  $(s+j)$-th Bernstein polynomial  of the fresco $\mathcal{F}_{f, \omega}^{[\alpha]}$ for some $j \in \mathbb{N}$.\\
Moreover, if $\xi_s = \xi_{s+1} = \cdots = \xi_{s+p}$ then there exists at least $p$ distinct values of $j \in \mathbb{N}$ such that $\xi_s$ is root of the $(s+j)$-th Bernstein polynomial of the fresco $\mathcal{F}_{f, \omega}^{[\alpha]}$.
\end{cor}

\parag{Proof} The proof of the first assertion is analogous to the proof of the theorem above. \\
The second assertion is an immediate consequence of the fact that the roots of the Bernstein polynomial of a semi-simple fresco are simple, applied to the successive semi-simple quotients
$$ S_d(\mathcal{F}^{[\alpha]})\big/S_{d-1}(\mathcal{F}^{[\alpha]}) $$
for $d =  s+1,  s+2, \dots, s+p$.$\hfill \blacksquare$\\

Now we conclude by a theorem  which combines the previous  results to precise the link between the first pole of order $\geq p$in $-\alpha - \mathbb{N}$ for a given pair $(\omega, \omega')$ with the roots of the Bernstein polynomials of order $\geq p$ of the fresco $\mathcal{F}_{f,\omega}$ associated to $(f,\omega)$.

\begin{thm}\label{sauts}
In the standart situation, assume that the hypothesis $H(\alpha, 1)$ is satisfied. Let $\omega$ be in $\Omega^{n+1}_0$ and define the fresco
 $\mathcal{F}_\omega := B[a] \omega$. Assume that  $p := d(\mathcal{F}_\omega^{[\alpha]})$ is at least equal to $1$  and choose\footnote{such a $\omega'$ exists thanks to Theorem \ref{Higher Bernstein 1}.} $\omega' \in \Omega^{n+1}_0$ such that there exists $h \in \mathbb{Z}$ such that $F^{\omega, \omega'}_h(\lambda)$ has a pole of order $p$ at some point in $-\alpha - \mathbb{N}$. For each $j \in [1, p]$, let $m_j$ be the integer which has the property $P(\omega, \omega', j)$. Then $-\alpha -m_j$ is a root of at least one of the polynomials $B^{j+q}(\mathcal{F}_\omega^{[\alpha]})$, for some integer $q$ in $\mathbb{N}$.
\end{thm}

\parag{Proof} Assuming that for some $j \geq 1$ no root of  the polynomials $B^{j+q}_{\F^{[\alpha]}}$ for $q \geq 0$ is equal to $-\alpha -m_j$ allows to find a J-H. sequence of $\mathcal{F}_\omega$  such that  the corresponding generator of the annihilator of $\omega$ is of the form $\Pi := \Pi_2\Pi_1$ where $\Pi_1$ has no factor $(a - \lambda_h b)$ with $\lambda_h + h - k$ equal to $\alpha + m_j$,  where the nilpotent order of the fresco
$({B[a]\big/ B[a] \Pi_2})^{[\alpha]})$ is at most  $ j-1$ and where $k$ is the rank of $\mathcal{F}_\omega$. Then we conclude as in the previous Theorem using Corollary \ref{pas de pole}.$\hfill \blacksquare$\\

\section{ Examples}

 It is, in general, rather difficult to compute the Bernstein of the fresco associated to a given pair $(f, \omega)$, even in the case where $f$ has an isolated singularity. \\
 Nevertheless, in the case where $f$ is a polynomial in $\mathbb{C}[x_0, \dots, x_n]$ having $(n+2)$ monomials, we describe  in the article  \cite{[B.22]},  a rather elementary method to obtain an estimation for the Bernstein polynomial of the fresco $\mathcal{F}_{f,\, \omega}$ associated to a monomial $(n+1)$-form $\omega$.\\
Of course, when the full Bernstein polynomial has a  root of multiplicity $k \geq 2$ then this root is also a root of the $j$-th Bernstein polynomial for at least $k$ values of $j$ but when the Bernstein polynomial has only simple roots, the computation of the higher  Bernstein polynomials, even in the special situation of   \cite{[B.22]}, is not easy. We present below some examples where the second Bernstein polynomial  is not trivial but where the full  Bernstein polynomial has no multiple root.

\begin{prop}\label{example 0}
Let $f(x,y,z) := xy^3 + yz^3 + zx^3 + \lambda xyz $ where $\lambda \not= 0$ is any complex number  which is a parameter, and consider the holomorphic forms
 $$\omega_1 := dx\wedge dy\wedge dz,\  \omega_2 := y^3z^2 \omega_1, \ \omega_3 =  y^7 \omega_1, \ {\rm and}  \   \omega_4 := xy^3 \omega_1.$$
 Then, in each of these cases, the fresco $\mathcal{F}_{f, \,\omega_i}$ is a rank $2$ theme and the second Bernstein polynomial  is equal  respectively to $x+1$, $x+4$, $x+5$ and $x+3$.\\
  Moreover, for $i = 3, 4$ the corresponding (full) Bernstein polynomial of the corresponding frescos  has only simple roots.
 \end{prop}

Note that this proposition allows to apply Theorem \ref{improve 1} to conclude that for each $i \in \{1,2,3,4\}$, there exists some integer $h$ and some germ $\omega_i' \in \Omega^3_0$ such that  the meromorphic extension of 
$$ F_h^{\omega_i,\omega_i'}(\lambda) = \frac{1}{\Gamma(\lambda)}\int_X \vert f\vert^{2\lambda} \bar f^{-h} \rho \omega_i\wedge \bar\omega_i' $$
has a double pole at  the point $ \lambda_i$ equal to the root of the second Bernstein polynomial of the fresco $\mathcal{F}_{f,\omega_i}$.\\

The proof of this proposition uses several lemmas and the technic of computation described in \cite{[B.22]} (see paragraph 4.3.2 in {\it loc. cit.}).

\begin{lemma}\label{gen. theme}
Let $e$ be a generator of the rank $2$ theme $T := B[a]/B[a](a-2b)(a-b)$ (which is the unique fresco  with Bernstein polynomial $(x+1)^2$). Assume that we have three  homogeneous polynomials 
$P, Q$ and $R$ in $B[a]$ of  respective degrees $3,4$ and $k$ with the following conditions
\begin{enumerate}
\item  $P,Q $ and $R$ are monic in $a$.
\item Then exists a non zero constant $c$ such that $P + cQ$ kills $e$ in $T$.
\item The Bernstein polynomial of $Q$\footnote{By definition $B_P$ is defined by the formula $$ (-b)^pB_P(-b^{-1}a) = P$$ where $P$ is in $B[a]$, is homogeneous in (a,b) of degree $p$ and monic in $a$. This is the Bernstein polynomial of the fresco  $B[a]/B[a]P$.} is not a multiple of $(x+1)$ or of  $(x+2)$.
\item The Bernstein polynomial of $R$ is not a multiple of $(x+3)(x+2)(x+1)$
\end{enumerate}
Then $Re$ generates a rank two sub-theme in $T$.
\end{lemma}

\parag{Proof} First, remark that our hypothesis implies that $P = (a - \nu b)(a - 2b)(a - b)$ for some $\nu \in \mathbb{C}$ since $T$ is isomorphic to $B[a]\big/B[a](a-2b)(a-b) $. We may realize $T$ in the simple pole asymptotic expansion module with rank $2$ which is isomorphic to $T^\sharp$
$$S\Xi_1^{(1)} = \Xi_1^{(2)}\big/\Xi_1^{(0)} \simeq \mathbb{C}[[s]](Log\, s)^2 \oplus \mathbb{C}[[s]Log\, s$$
 where $a$ is the multiplication by $s$ and $b$ is defined by  $ab - ba = b^2$ and 
$$ b(Log\, s)= sLog\, s \quad {\rm and} \quad b((Log\,s)^2) = s(Log\, s)^2 - 2sLog\, s .$$ 
Then let us prove that  image of $e$ in $\Xi_1^{(2)}\big/\Xi_1^{(0)} $ may be written 
\begin{equation*}
 e = u(Log\, s)^2 + vs(Log\, s)^2 + ws^3(Log\, s)^2 +  s^4\mathbb{C}[[s]] (Log\, s)^2+ \mathbb{C}[[s]](Log\, s) \tag{@}
\end{equation*}
where $\varphi$ is in $\Xi_1^{(2)}\big/\Xi_1^{(0)} $ and where $uvw\not= 0$ are complex numbers. \\
  Remark that  the only restrictive condition for writing $e$ as in $(@)$ is the condition  $uvw\not= 0$. The condition $u \not= 0$ is easy because we assume that $e$ is a generator of $T$ with Bernstein polynomial $(x+1)^2$, so writing $e$ as a $\mathbb{C}[[b]]$-linear combination of the $\mathbb{C}[[b]]$-basis $e_1 = (Log\, s)^2$ and $e_2 = Log\, s$ of $T$  we see that the coefficient of $e_1$ must be invertible in $\mathbb{C}[[b]]$.\\
 But the condition $(P + cQ)(e) = 0$ implies, since the Bernstein element of $T$ is $(a-2b)(a-b)$, that we  may write\footnote{In our choice of $f$ and $\omega_1$, $\nu = 3$.} $P = (a-\nu b)(a - 2b)(a - b)$. \\
 The annihilator of  $(Log\, s)^2$ in $\Xi_1^{(2)}\big/\Xi_1^{(0)} $ is the ideal $B[a](a-2b)(a-b)$ so we have $P( (Log\, s)^2 ) = 0$ in $T$. Since  $Q((Log\, s)^2)$ has a non zero term in $s^4(Log\, s)^2$, because $-1$ is not a root of $B_Q$, only the term coming from
 $$ P(s(Log\, s)^2) = \frac{4-\nu}{24}s^4(Log\, s)^2 \quad {\rm modulo} \ \mathbb{C}[[s]]Log\, s $$
 can compensate for this term, in order  to obtain  the equality $(P + cQ)(e) = 0$. Then $u\not= 0$. implies  $v \not= 0$. \\
 But now, the only term which can kill the non zero term in $s^5(Log\, s)^2$ coming from $Q(vs(Log\, s)^2)$ (using  that $B_{Q}$ is not a multiple of $(x+2)$)  can only come from $P(ws^2(Log\,s)^2)$ and this proves that $w \not= 0$. So the assertion  $(@)$ holds true.\\
 Now if $R$ is homogeneous of degree $k$ in $(a,b)$ a necessary condition on $R$ such that $R(e)$ has no term in $s^{k+i}(Log\, s)^2$,  for $i= 0,1,2$,  is that $B_R$ divides $(x+1)(x+2)(x+3)$. So, when it is not the case, Lemma  5.2.4 in \cite{[part I]}  implies that $R(e)$ is a rank $2$ theme and that its second Bernstein polynomial has a (unique)  root equal to $-(k+j)$ where $-j$ is the smallest integer among  $\{-1,-2,-3\}$ which is not a root of $B_R$. $\hfill \blacksquare$\\
 
 Note that the Lemma above may be easily generalized  to many $[\alpha]$-primitive frescos  provided that the nilpotent order is known and that it has  a  generator which admits a enough simple element in $B[a]$ belonging to its annihilator.\\
   
  \begin{lemma}\label{example 1}
 In the situation of Proposition \ref{example 0}, the frescos generated by the forms 
 $$\omega_1 := dx\wedge dy\wedge dz,  \ \omega_2 := y^3z^2 \omega_1, \ \omega_3 := y^7\omega_1, \  {\rm and} \  \omega_4 := xy^3\omega_1 $$ 
 generate rank $2$ $[1]$-primitive themes. Their Bernstein polynomials are respectively equal to $$(x+1)^2, \quad (x+3)^2 \ {\rm or} \ (x+2)(x+3), \quad  (x+3)(x+5) \quad {\rm and} \quad (x+2)(x+3)$$
 and  their respective $2$-Bernstein polynomials are  $(x+1)$, $(x+3)$, $(x+5)$ and $(x+3)$. In the  cases $i = 3,4$  there is no double root for the Bernstein polynomial of $\mathcal{F}_{f, \omega_i}.$ 
 \end{lemma}

 \parag{Proof} The first point is to show that $\mathcal{F}_{f,\omega_1}$ has  rank $2$. Since $f$ has an isolated singularity at the origin, we have $Ker df^{n}= df\wedge \Omega^{n-1}$ and then 
 $H_0^{n+1}/bH_0^{n+1} \simeq \mathcal{O}_{0}/J(f)$ and $H_0^{n+1}$ has no $b$-torsion and no $a$-torsion. Since $f$ is not\footnote{This point is not so easy to check directly. But  the rank is not $1$ since this would implies that  this fresco  has a simple pole and the argument used in  Lemma  \ref{gen. theme} gives then a contradiction.} 
   in $J(f)$ the image of $\omega_1$ and $a\omega_1 = f\omega_1$ in $H^{n+1}$ are linearly independent (over $\mathbb{C}$) and then the rank of $B[a]\omega_1$ is at least equal to $2$. Now the computation in \cite{[B.22]} (see 4.3.2) shows that the Bernstein polynomial of this fresco divides $(x+1)^3$ (see also the detailed computation below). So it is a theme of rank $2$ or $3$. But using our main result, the rank $3$ would imply that there exists a pole of order $3$ for some $F^{\omega_1,\omega'}_h(\lambda)$ which is impossible\footnote{This would give an order $4$ pole for the meromorphic continuation of $\vert f\vert ^{2\lambda}$ !} in $\mathbb{C}^3$. So $\mathcal{F}_{f,\omega_1}$ is a rank $2$ theme with Bernstein polynomial $(x+1)^2$. The computation in \cite{[B.22]} gives that $P_3 +  c\lambda^{-4}P_4$ kills $\omega_1$ in $H_0^{n+1}$ where 
 $$ P_3 := (a-3b)(a-2b)(a-b), P_4 = (a -(13/4)b)(a-(5/2)b)(a-(7/4)b)a, {\rm and } \ c = 4^4 $$
 This is easily obtained by using the technic of the computation of {\it loc.cit.} (see  the  detailed computation in the Appendix of \cite{[B.23]}). 
 Then we may apply  Lemma \ref{gen. theme} to see that $\lambda m_1m_2 \omega_1 = \lambda (a-2b)(a-b)\omega_1$  generates rank $2$ themes in $H_0^{n+1}$.
 But the identity $\lambda m_1m_2 = m_4 y^3z^2$ shows that $\omega_2$ generates also  rank $2$ in $H_0^{n+1}$ since $m_4\omega_2 = \lambda m_1m_2\omega_1 = \lambda(a-2b)(a-b)\omega_1$ applying Lemma  \ref{gen. theme} with $R = (a-2b)(a-b)$ whose Bernstein polynomial is $(x+1)^2$. Moreover we see that $Re$ has a non zero term in $s^3(Log\, s)^2$.\\
 Since $m_4\omega_2$ generates a rank $2$ theme, then $\omega_2$ generates a rank $2$ theme also (the rank $3$ is again excluded because  it would imply that $f^2 \not\in J(f)$ which is impossible as explained above).\\
 The technic of computation in \cite{[B.22]}  applied to $\omega_2$ gives now that the Bernstein polynomial of the rank $2$ theme $B[a]\omega_2$ has to divide\footnote{This computation gives that $Q_3 + d\lambda^{-4}Q_4$ kills $\omega_2$ in $H_0^{n+1}$ with $Q_3 := (a-4b)(a-4b)(a-3b)$.} the polynomial $(x+2)(x+3)^2$.\\
 But the fact that $m_4\omega_2$ has a  non zero term in $s^3(Log\, s)^2$ (and no term in $(Log\, s)^2$ or in $s(Log\,s)^2$)  implies, since we have 
  $$m_4\omega_2 = 4(a - 2b)\omega_2$$
  $\omega_2$ has a non zero term in  $s^2(Log\, s)^2$  and then $-3$ is a root of the second Bernstein polynomial of the fresco $\mathcal{F}_{f, \omega_2}$. So the Bernstein polynomial is either $(x+2)(x+3)$ or $(x + 3)^2$.\\
   We know that   the Bernstein polynomial of $\mathcal{F}_{f,\omega_3}$ divides $(x+5)(x+3)(x+2)$ by using the technic of  \cite{[B.22]}. \\
   We know also  that $m_1^2m_4\omega_1 = \lambda m_3\omega_3$ has a non zero term in  $s^5(Log\, s)^2$ (as a consequence of Lemma \ref{gen. theme}) and, since  $-m_3\omega_3 = (a - 2b)\omega_3$ implies that $\omega_3$ has a non zero term in $s^4(Log\, s)^2$, the second Bernstein polynomial of $\mathcal{F}_{f, \omega_3}$ is $x+5$.\\
   
  Note that the Bernstein polynomial of the fresco $\mathcal{F}_{f, \omega_3}$ has two simple roots.\\
  
  The last  case is similar, since we know that $m_1\omega_1$  has a non zero term in $s^2(Log\, s)^2$.  So our assertion is consequence of the estimation of the Bernstein polynomial.
   $\hfill \blacksquare$\\
   
   The reader may find more details on the previous computations in \cite{[B.23]}.\\

\bigskip

{\it The author has no relevant financial or non-financial interests to disclose.
The author has no conflicts of interest to declare that are relevant to the content of this article.
The author certifies that he has no affiliations with or involvement in any organization or entity with any financial interest or non-financial interest in the subject matter or materials discussed in this manuscript.
The author has no financial or proprietary interests in any material discussed in this article}.

\end{document}